\documentclass[11pt,a4paper]{amsart}
\usepackage{amsfonts,layout}

\newcount\minute
\newcount\hour
\def\currenttime{%
	\minute\time
	\hour\minute
	\divide\hour60
	\the\hour:\multiply\hour60\advance\minute-\hour\the\minute}
\def\draftnote{{\it \today \quad  \currenttime \hfill  tex-file :   \jobname}}

\catcode`@=11
\@addtoreset{equation}{section}

\catcode`@=12
\newtheorem{Theorem}{Theorem}[section]
\newtheorem{Definition}{Definition}[section]
\newtheorem{Proposition}{Proposition}[section]
\newtheorem{Lemma}{Lemma}[section]

\newtheorem{Remark}{Remark}[section]

\newtheorem{Hyp.}{Hyp.}[section]

\begin{document}
%\draftnote  \bigskip \bigskip

\title[]{The cost of controlling degenerate parabolic equations by boundary controls}

\author{P. Cannarsa} 
\address{Dipartimento di Matematica, Universit\`a di Roma "Tor Vergata",
Via della Ricerca Scientifica, 00133 Roma, Italy}
\email{cannarsa@mat.uniroma2.it}

\author{P. Martinez} 
\address{Institut de Math\'ematiques de Toulouse, UMR CNRS 5219, Universit\'e Paul Sabatier Toulouse III \\ 118 route de Narbonne, 31 062 Toulouse Cedex 4, France} \email{Patrick.Martinez@math.univ-toulouse.fr}

\author{J. Vancostenoble} 
\address{Institut de Math\'ematiques de Toulouse, UMR CNRS 5219, Universit\'e Paul Sabatier Toulouse III \\ 118 route de Narbonne, 31 062 Toulouse Cedex 4, France} \email{Judith.Vancostenoble@math.univ-toulouse.fr}

%\date{\today}
\thanks{This research was partly supported by the Institut Mathematique de Toulouse and Istituto Nazionale di Alta Matematica through funds provided by the national group GNAMPA and the GDRE CONEDP. Moreover, this work was completed while the first author was visiting the Institut Henri Poincar\'e and Institut des Hautes \'Etudes Scientifiques on a CARMIN senior position.}

%%%%%%%%%%%%%%%%%%%%%%%%%%%%%%%%%%%%%%%%
%%%%%%%%%%%%%%%%%%%%%%%%%%%%%%%%%%%%%%%%
\begin{abstract}
We consider the one-dimensional degenerate parabolic equation 
$$ u_t - (x^\alpha u_x)_x =0 \qquad x\in(0,1),\ t \in (0,T) ,$$
controlled by a boundary force acting at the degeneracy point $x=0$.

First we study the reachable targets at some given time $T$
using $H^1$ controls, extending the moment method developed by Fattorini and Russell \cite{FR1, FR2} to this class of degenerate equations.

Then we investigate the controllability cost to drive an initial condition to rest, deriving optimal bounds with respect to $\alpha$ and deducing that the cost blows up as $\alpha \to 1^-$.

\end{abstract}
\maketitle

%\tableofcontents
%%%%%%%%%%%%%%%%%%%%%%%%%%%%%%%%%%%%%%%%
%%%%%%%%%%%%%%%%%%%%%%%%%%%%%%%%%%%%%%%%
\section{Introduction}

Null controllability of nondegenerate parabolic equations is by now well understood, either by locally distributed control or by a control acting on a part of the boundary, and we refer the reader to the seminal papers of Fattorini and Russell \cite{FR1, FR2}, and of Fursikov and Imanuvilov 
\cite{Fursikov}. 

However, many problems that are relevant for applications are described by degenerate equations, with degeneracy occurring at the boundary of the space domain. We can mention the question of invariant sets for diffusion process in probability, the study of the velocity field of a laminar flow on a flat plate, the Budyko-Sellers climate models, the Fleming-Viot gene frequency model, referring, e.g., to 
\cite{memoire} for details.
The typical example of such degenerate problems is the equation
\begin{equation}
\label{eq-deg}
u_t - (x^\alpha u_x)_x =0 \qquad x\in(0,1),\ t>0,
\end{equation}
where $\alpha >0$ is given. The diffusion coefficient in \eqref{eq-deg} vanishes at $x=0$.
The question of null controllability by a locally distributed control has been solved in \cite{sicon2008}, and developped in several directions (more general classes of degenerate equations (\cite{fatiha, non-div0, non-div, Martinez}), and in space dimension 2 (\cite{memoire})). 

The question of null controllability by a boundary control acting where the degeneracy occurs, has been studied recently:
\begin{itemize}
\item Cannarsa-Tort-Yamamoto \cite{CTY} proved an approximate controllability result;
\item Gueye \cite{Gueye} proved the expected null controllability result, studying the degenerate wave equation and then applying the transmutation method (\cite{Ervedoza1, Ervedoza2});
\item Martin, Rosier and Rouchon \cite{Rosier} proved a null controllability result using the flatness approach, which can be applied to general situations (degenerate or singular parabolic equations), and gives the control and the solution as series.
\end{itemize}

The main goal of this paper is to study the dependence of the controllability properties with respect to the degeneracy parameter $\alpha$, in the typical example \eqref{eq-deg} when the control acts on the boundary at the degeneracy point $x=0$:

\begin{itemize} 

\item our first result concerns the reachable targets using $H^1$ controls: we prove that there is an explicit subset $\mathcal{P}_{\alpha,  T} \subset L^2(0,1)$, which is dense in $L^2(0,1)$, such that every $u_T \in \mathcal{P}_{\alpha,  T}$ is reachable with $H^1$ controls (see Theorem \ref{thm2});

\item since the reachable set $\mathcal{R}_{\alpha,  T}$ contains a subset $\mathcal{P}_{\alpha,  T}$ which is dense in $L^2(0,1)$, it is interesting to look for targets that could be reached for all parameter $\alpha \in [0,1)$; 
however, we prove that $\bigcap _{\alpha \in [0,1)} \mathcal{P}_{\alpha,  T} = \{0 \}$, hence $0$ is the only target that we are sure that can be reached for all parameter $\alpha \in [0,1)$ (see Proposition \ref{prop-always-targets});

\item since $0$ is reachable for all parameter $\alpha \in [0,1)$, it is interesting to measure the cost to drive an initial condition to $0$ with respect to the parameter $\alpha$;
we prove that the controllability cost blows up with order $\frac{1}{1-\alpha}$ as the degeneracy parameter goes to $1^-$, providing (optimal) upper and lower bounds (see Theorem \ref{thm-cost-bd}).

\end{itemize}

The proofs are based on the moment method developed by Fattorini and Russel \cite{FR1, FR2}. We extend their method and some of their results to this degenerate case, and then we take advantage of the explicit expressions of the control that it provides (in terms of Bessel functions and their zeros) to obtain upper and lower bounds of the null controllability cost.

\medskip

The paper is organized as follows: in section \ref{s-results} we state precisely our results; in section \ref{s3}, we summarize the definitions and properties of Bessel functions that are useful to solve the Sturm-Liouville problem;  section \ref{s7} is devoted to the proof of Theorem \ref{thm2} (concerning the subpart $\mathcal{P}_{\alpha,  T}$ of reachable targets); section \ref{sec-struct} is devoted to the proof of Proposition \ref{prop-always-targets} (concerning $\bigcap _{\alpha \in [0,1)} \mathcal{P}_{\alpha,  T}$); section \ref{sec-cost} is devoted to the proof of Theorem \ref{thm-cost-bd} (concerning the cost of null controllability).

%%%%%%%%%%%%%%%%%%%%%%

\section{Setting of the problem and main results}
\label{s-results}

%%%%%%%%%%%%%%%%%%%

We are interested in the controllability properties of the problem

\begin{equation}
\label{eq-u-G}
\begin{cases}
u_t - (x^\alpha u_x)_x = 0, \\
u(0,t)=G(t), \\
u(1,t) = 0, \\
u(x,0)=u_0(x) ,
\end{cases}
\end{equation}
that is when the control acts at the degeneracy point $0$ through a nonhomogeneous Dirichlet boundary condition.
First, we recall that, in a general way, the well-posedness of degenerate parabolic equations is stated in weighted Sobolev spaces. We will consider the problem when $\alpha \in [0,1)$; in this case, the Dirichlet boundary control makes sense, as we explain in the following.

\subsection{Well-posedness and eigenvalue problem}

%%%%%%%%%%%%%%

\subsubsection{A preliminary transformation}

To define the solution of the boundary value problem \eqref{eq-u-G}, we transform it into a problem with homogeneous Dirichlet boundary conditions and a source term (depending on the control $G$): 
consider 
$$ p (x) := \int _x ^1 \frac{ds}{s^\alpha} .$$
Then, formally, if $u$ is a solution of \eqref{eq-u-G}, then the function $v$ defined by
\begin{equation}
\label{v:=u,G}
 v(x,t) = u(x,t) - \frac{p(x)}{p(0)} G(t)  = u(x,t) -  (1-x^{1-\alpha}) G(t) 
\end{equation}
satisfies the auxiliary problem
\begin{equation}
\label{eq-v-G'}
\begin{cases}
v_t - (x^\alpha v_x)_x = - \frac{p(x)}{p(0)} G'(t) , \\
v(0,t)=0, \\
v(1,t) = 0, \\
v(x,0)=u_0(x) - \frac{p(x)}{p(0)} G(0).
\end{cases}
\end{equation}
Reciprocally, given $g \in L^2(0,1)$, consider the solution $v$ of 
\begin{equation}
\label{eq-v-g}
\begin{cases}
v_t - (x^\alpha v_x)_x = - \frac{p(x)}{p(0)} g(t) , \\
v(0,t)=0, \\
v(1,t) = 0, \\
v(x,0)=v_0(x) .
\end{cases}
\end{equation}
Then the function $u$ defined by
\begin{equation}
\label{u:=v,g} 
u(x,t) = v(x,t) + \frac{p(x)}{p(0)} \int _0 ^t g(\tau) \, d\tau  
\end{equation}
satisfies
\begin{equation}
\label{eq-u-g}
\begin{cases}
u_t - (au_x)_x = 0, \\
u(0,t) = \int _0 ^t g(\tau) \, d\tau , \\
u(1,t) = 0, \\
u(x,0)=v_0(x) .
\end{cases}
\end{equation}
This motivates the following definition of what is the solution of the boundary value problem \eqref{eq-u-G}, as we explain in the following.

%%%%%%%%%%%%%%

\subsubsection{Well-posedness of the problem for $H^1(0,T)$ boundary controls}

Consider $\alpha \in [0,1)$ and
$$ H^1 _\alpha (0,1) := \{ u\in L^2(0,1), \text{ $u$ absolutely continuous on $[0,1]$}, x^{\alpha /2} u' \in L^2(0,1) \}, $$
$$ H^1 _{\alpha,0} (0,1) := \{ u \in H^1 _\alpha (0,1), u(0)=0=u(1) \} ,$$
$$ H^2 _\alpha (0,1) := \{u \in H^1 _\alpha (0,1), x^\alpha u' \in H^1 (0,1) \} ,$$
and the unbounded operator $A:D(A)\subset L^2(0,1)\to L^2(0,1)$ defined by
\begin{equation*}
%\label{D(A)-1}
\begin{cases}
\forall u \in D(A), \quad  Au:= (x^\alpha u_x)_x,  & \\
D(A) :=\{ u \in H^1_{\alpha ,0}(0,1)  \ \mid \ x^\alpha u_x \in  H^1(0,1) \}.&
\end{cases}
\end{equation*}
In both cases, the following results hold, (see, e.g., \cite{campiti} and \cite{CMV4}).

\begin{Proposition}
\label{prop-A} $A: D(A) \subset L^2(0,1)\to L^2(0,1)$ is a self-adjoint negative operator
with dense domain.
\end{Proposition}
Hence, $A$ is the infinitesimal generator of an analytic semigroup of contractions $e^{tA}$
on $L^2(0,1)$. 
Given a source term $h$ in  $L^2((0,1)\times (0,T))$ and an initial condition $v_0 \in  L^2(0,1)$, consider the problem
\begin{equation}
\label{eq-v-h}
\begin{cases}
v_t - (x^\alpha v_x)_x = h(x,t), \\
v(0,t)=0, \\
v(1,t) = 0, \\
v(x,0)=v_0(x)  .
\end{cases}
\end{equation}
 The function $v \in \mathcal C ^0 ([0,T]; L^2(0,1))  \cap L^2(0,T; H^1_{\alpha,0} (0,1))$
given by the variation of constant formula
$$ v(\cdot ,t) = e^{tA} v_0 + \int _0 ^t e^{(t-s)A} h(\cdot, s) \, ds
$$ 
is called the mild solution of \eqref{eq-v-h}. We say that a function 
$$v \in
\mathcal C ^0 ([0,T]; H^1_{\alpha,0} (0,1))  \cap   H^1(0,T; L^2(0,1)) \cap    L^2 (0,T; D(A)) $$
is a strict solution of \eqref{eq-v-h} if $v$ satisfies $v_t - (x^\alpha v_x)_x = h(x,t)$ almost everywhere
in $(0,1)\times (0,T)$, and the initial and boundary conditions  for all $t\in [0,T]$ and all $x\in [0,1]$.

\begin{Proposition}
If $v_0 \in H^1 _{\alpha,0}(0,1)$, then the mild solution of \eqref{eq-v-h} is the unique strict solution of \eqref{eq-v-h}.
\end{Proposition}

In particular, the above notions of solution apply to problem \eqref{eq-v-g} taking $h(x,t) = -\frac{p(x)}{p(0)} g(t)$ and $v_0(x) = u_0(x)-\frac{p(x)}{p(0)}G(0)$ . This allows us to define in a suitable way the solution $u$ of \eqref{eq-u-G}:

\begin{Definition}
a) We say that $u \in C ([0,T]; L^2(0,1)) \cap L^2(0,T; H^1 _\alpha (0,1))$ is the mild solution of 
\eqref{eq-u-G} if $v$ defined by \eqref{v:=u,G} is the mild solution of \eqref{eq-v-G'}.

b) We say that $$u \in C ([0,T]; H^1 _{\alpha} (0,1)) \cap H^1(0,T; L^2 (0,1))\cap L^2(0,T; H^2 _\alpha (0,1))$$ is the strict solution of \eqref{eq-u-G} if $v$ defined by \eqref{v:=u,G} is the strict solution of \eqref{eq-v-G'}.
\end{Definition}
Then we immediately obtain
\begin{Proposition}
\label{prop-u-G}
a) Given $u_0 \in L^2 (0,1)$, $G\in H^1(0,T)$, problem \eqref{eq-u-G} admits a unique mild solution.

b) Given $u_0 \in H^1 _{\alpha} (0,1)$, $G \in H^1(0,T)$ such that $G(0)=u_0(0)$, problem \eqref{eq-u-G} admits a unique strict solution. In particular it holds true when
$u_0 \in H^1 _{\alpha,0}(0,1)$, $G \in H^1(0,T)$ and $G(0)=0$.
\end{Proposition}
\noindent The proof of Proposition \ref{prop-u-G} follows immediately, noting that
$$ \tilde G (x,t) := \frac{p(x)}{p(0)} G(t) $$ satisfies
$$ \tilde G \in C ([0,T]; H^1 _{\alpha} (0,1)) \cap H^1(0,T; L^2 (0,1))\cap L^2(0,T; H^2 _\alpha (0,1)) .$$

\begin{Remark}\label{rq-reg} {\rm 
When $u_0 \in L^2(0,1)$, given $\tau >0$, we have $v(\cdot, \tau)\in  H^2 _\alpha (0,1) \cap H^1 _\alpha (0,1)$, therefore the solution of \eqref{eq-v-g} is strict on $[\tau,T]$. The same is true for the solution of \eqref{eq-u-G}.} 
\end{Remark}

%%%%%%%%%%%%%%

\subsubsection{The eigenvalue problem and the associated eigenfunctions}

The knowledge of the eigenvalues and associated eigenfunctions of the degenerate diffusion operator $ y \mapsto - (x^\alpha y')'$, i.e. the solutions $(\lambda, y)$ of
\begin{equation}\label{pbm-vp}
\begin{cases}
 - (x^\alpha y'(x))' =\lambda y(x) & \qquad x\in (0,1),\\
y (0)=0 ,\\ 
y(1)=0. & 
\end{cases}
\end{equation}
 will be essential for our purposes.
It is well-known that Bessel functions play an important role in this problem, see, e.g., Kamke \cite{Kamke}. For $\alpha \in [0,1)$, let
$$ \nu _\alpha := \frac{ 1-\alpha }{2-\alpha}, 
\qquad \kappa _\alpha:= \frac{2-\alpha}{2}.$$
Given $\nu \geq 0$, we denote by $J_\nu$ the Bessel function of first kind and of order $\nu$ (see section \ref {s3}) and denote $j_{\nu,1}< j_{\nu,2} < \dots < j_{\nu,n} <\dots$ the sequence of positive zeros of $J_\nu$. 
Then the admissible eigenvalues $\lambda$ for problem \eqref{pbm-vp} are given by 
\begin{equation}
\label{vp}
\forall n \geq 1, \qquad \lambda_{\alpha, n} =  \kappa _\alpha ^2 j_{\nu _\alpha ,n}^2
\end{equation}
and the corresponding normalized (in  $L^2(0,1)$) eigenfunctions takes the form
\begin{equation}
\label{fp}
 \Phi_{\alpha, n}(x)=  \frac{\sqrt{2 \kappa _\alpha }}{\vert J'_{\nu_\alpha} (j_{\nu_\alpha,n} ) \vert} 
x^{(1-\alpha)/2} J_{\nu _\alpha} (j_{\nu_\alpha,n} x ^{\kappa_\alpha}), \qquad x \in (0,1).
\end{equation}
Moreover the family $(\Phi_{\alpha, n})_{n\geq 1}$ forms an orthonormal basis of $L^2(0,1)$.

\begin{Remark}\label{rq1} 
{\rm 
Let us observe that in the case $\alpha =0$ (which corresponds to the non degenerate heat equation), we have $\nu _0=1/2$ and $\kappa _0=1$. Classical properties of the Gamma function (\cite[section 1.2, formula (1.2.3), p. 3]{Lebedev}) show that 
$$J_{1/2} (x)= \frac{\sqrt{2}}{\sqrt{\pi x}} \sin x ;$$
in that case, we deduce from \eqref{vp} that $ \lambda_{0,n}=(n\pi)^2$ and from \eqref{fp} that $\Phi_{0,n}(x)=\sqrt{2} \sin(n\pi x)$ for all $n \geq 1$ and $x\in (0,1)$. Therefore one recovers the well-known eigenvalues and eigenfunctions of the Laplace operator in $(0,1)$.
}
\end{Remark}

%%%%%%%%%%%%%%

\subsection{Main results: the reachable set and the cost of null controllability}
\label{sub-resultats}

\subsubsection{The controllability problem}

The first problem we address concerns the boundary controllability of equation \eqref{eq-u-G} using a control acting at the degeneracy point. Given $\alpha \in [0,1)$, $T>0$, $u_0, u_T \in L^2(0,1)$, we wish to find $G \in H^1 ( 0,T)$ that drives the solution $u$ of \eqref{eq-u-G} from $u_0$ to $u_T$ in time $T$. (Of course, due to the regularizing effect, it is clear that we will not be able to reach targets $u_T$ with low regularity.)

We will also be interested in the targets that can be reached for all the parameters $\alpha$. If such a target exists, it makes sense to evaluate the cost to reach it with respect to the degeracy parameter $\alpha$.

%%%%%%%%%%

\subsubsection{A Fourier-Bessel description of the reachable set}

The result we obtain is the following one: 

\begin{Theorem}\label{thm2}
Let $\alpha \in [0,1)$. Consider a target function $u_T$, and the sequence $(\mu_{\alpha ,n} ^T)_{n\geq 1}$ of  its Fourier coefficients:
$$\mu_{\alpha ,n} ^T = \int _0 ^1 u_T (x) \Phi _{\alpha,n} (x) \,dx .$$
Then there exists some $K>0$ independent of $\alpha \in [0,1)$ such that, if
\begin{equation}
\label{decr-target-gene}
\sum_{n \geq 1} n^{3/2} \vert \mu _{\alpha ,n} ^T\vert  e^{K \kappa_\alpha \pi n} < \infty ,
\end{equation}
then $u_T$ is a reachable target: given
$T>0$ and $u_0 \in L^2(0,1)$, there exists $G_\alpha \in H^1 ((0,1)\times (0,T))$ such that the solution $u$ of \eqref{eq-u-G} controlled by $G_\alpha$ satisfies $u(T) =u_T.$
\end{Theorem}

We will denote by
$\mathcal{P}_{\alpha,  T}$ the set of $u_T$ that satisfy \eqref{decr-target-gene}.

\begin{Remark} \hfill {\rm

\begin{itemize}
\item Of course condition \eqref{decr-target-gene} is satisfied if $\mu_{\alpha ,n} ^T=0$ for all $n$ large enough. Since finite linear combinations of $\Phi _{\alpha, n}$ are dense in $L^2(0,1)$, the reachable targets form a dense subset of $L^2(0,1)$, which was already known from \cite{CTY}. Our result allows us to be more precise on the reachable targets. 

\item We underline the fact that \eqref{decr-target-gene} is independent of $T$. Indeed, it is well-known in a general setting that the reachable set $\mathcal R_T$ of the targets that can be attained at time $T$ does not depend on $T$, see Seidman \cite{Seidman}. 

\item A condition like \eqref{decr-target-gene} already appears in the pioneering works of Fattorini and Russell \cite{FR1, FR2}. Ervedoza and Zuazua \cite{Ervedoza1} proved a similar (in fact a slightly better) condition in a general context. One could provide an explicit estimate of the constant $K$ that appears in \eqref{decr-target-gene}, but we emphasize the fact that it does not depend on $\alpha$. This is worth to be noted since we are interested in the behavior of the reachable set with respect to the degeneracy parameter $\alpha$. 
\end{itemize} }

\end{Remark}

%%%%%%%%%%%%%%

\subsubsection{The regularity of the targets and the question of the targets that are reachable for all $\alpha \in [0,1)$}

Fattorini and Russell \cite{FR1} noted that in the case of the heat equation, i.e. when $\alpha =0$, a reachable target is the restriction to $[0,1]$ of an analytic function. 
Let us study what can be said in our case. We prove the following regularity result:

\begin{Proposition}
\label{prop-reg-targets}
Consider a sequence $(\mu_{\alpha ,n} ^T)_{n\geq 1}$ such that, for some $K>0$, the sequence $(\mu_{\alpha ,n} ^T e^{Kn})_{n\geq 1}$ is bounded.
Consider 
$$ \forall x \in [0,1], \quad u_T (x) := \sum _{n=1} ^\infty \mu_{\alpha ,n} ^T \Phi _{\alpha ,n} (x). $$
Then $u_T$ has the following property: there exists an even function $F_{\alpha}$, holomorphic in the strip $\{ z\in \Bbb C, \vert \Im z \vert < \frac{K}{\pi} \}$ such that 
\begin{equation}
\label{eq-structure-target}
\forall x\in [0,1], \quad u_T (x) = x^{1-\alpha} F_{\alpha} (x^{\kappa_\alpha}) .
\end{equation}
\end{Proposition}

This regularity result extends in a natural way the result of Fattorini and Russell \cite{FR1}, and it has the following consequences:

\begin{Proposition}
\label{prop-always-targets}

a) Given $\alpha \in [0,1)$, if $u_T \in \mathcal{P}_{\alpha,  T}$, i.e. if the Fourier coefficients of $u_T$ satisfy \eqref{decr-target-gene}, then $u_T$ is reachable, and there exists an even function $F_{\alpha }$, holomorphic in  the strip $\{ z \in \Bbb C, \vert \Im z \vert < \frac{K}{\pi} \}$ such that \eqref{eq-structure-target} holds.

b) The following property holds: 
$$\bigcap _{\alpha \in [0,1)} \mathcal{P}_{\alpha,  T} = \{0 \} .$$
 Hence, the only $u_T$ that satisfies \eqref{decr-target-gene} for all $\alpha \in [0,1)$ is $u_T=0$.
\end{Proposition}

\begin{Remark} {\rm The problem of establishing whether zero is the only target that can be reached for all $\alpha \in [0,1)$ is widely open.}
\end{Remark}

%%%%%%%%%%%%%%

\subsubsection{The cost of null controllability}

Finally, since $0$ can be reached for all $\alpha \in [0,1)$, it is interesting to measure the cost to drive any $u_0$ to $0$ in time $T$, with respect to $\alpha$.

We define the controllability costs in the following way:
given $u_0 \in L^2(0,1)$, we consider the set of admissible controls that drive the solution $u$ of \eqref{eq-u-G} to $0$ in time $T$:
$$ \mathcal U _{ad}(\alpha, u_0) := \{ G \in H^1(0,T), u^{(G)} (T)=0 \},$$
where $u^{(G)}$ denotes the solution of \eqref{eq-u-G}; we consider the controllability cost
\begin{equation}
\label{cost-fr-mm}
 C^{H^1} (\alpha, u_0) := \inf _{G\in \mathcal U _{ad}(\alpha, u_0)} \Vert G \Vert _{H^1(0,T)} ,
\end{equation}
which is the minimal value to drive $u_0$ to $0$.
We also consider a global notion of controllability cost:
\begin{equation}
\label{cost-fr}
 C^{H^1}_{bd-ctr}(\alpha) := \sup _{\Vert u_0 \Vert =1} C^{H^1} (\alpha, u_0).
\end{equation}
Similar notions were already being considered, see in particular Fernandez-Cara and Zuazua \cite{FCEZ}. 

Then we prove the following

\begin{Theorem}
\label{thm-cost-bd}

a) Given $u_0 \in L^2(0,1)$, there exists $M_1(u_0)$ independent of $\alpha \in [0,1)$, and $M_2$ independent of $u_0$ and $\alpha$ such that
\begin{equation}
\label{eq-cost-fr-u_0}
\frac{M_1(u_0)}{1-\alpha}  \leq C^{H^1} (\alpha, u_0) \leq \frac{M_2}{1-\alpha} \Vert u_0 \Vert _{L^2} .
\end{equation}

b) There exist two positive constants $M_1, M_2$ independent of $\alpha \in [0,1)$ such that
\begin{equation}
\label{eq-cost-fr}
\frac{M_1}{1-\alpha} \leq C^{H^1}_{bd-ctr}(\alpha) \leq \frac{M_2}{1-\alpha}.
\end{equation}
\end{Theorem}

\begin{Remark} {\rm
This shows that the controllability cost blows up as $\alpha \to 1^-$, and that our upper estimate is optimal.}
\end{Remark}

%%%%%%%%%%%%%%%%%

\subsection{Additional comments and related questions}

\subsubsection{The question of uniformly reachable targets} 
As in \cite{Ervedoza1}, and of course as in \cite{FR1}, we obtain a subpart $\mathcal P _{\alpha, T}$ of the reachable set  $\mathcal R _{\alpha, T}$, and we prove a somewhat negative result concerning $\bigcap _{\alpha \in [0,1)} \mathcal P _{\alpha, T}$ in Proposition \ref{prop-always-targets}. It would be interesting to obtain a result concerning $\bigcap _{\alpha \in [0,1)} \mathcal R _{\alpha, T}$.

%%%%%%%%%%%%

\subsubsection{The cost of null controllability}
It would be interesting to improve (if possible) \eqref {eq-cost-fr-u_0} of Theorem \ref{thm-cost-bd} to obtain a lower bound that depends on $\Vert u_0 \Vert _{L^2(0,1)}$.

%%%%%%%%%%%%

\subsubsection{The question of locally distributed controls}
In \cite{CMV-cost-loc}, we will study the same questions when the control is locally distributed in $(0,1)$.

%%%%%%%%%%%%%%%%%%%%%%%%%%%%%%%%%

\section{Useful tools from Bessel's theory for the Sturm-Liouville problem} \label{s3}

In this section, we recall existing tools, that we will need to prove our results, stated in section 
\ref{sub-resultats}. Note that  one can observe that if $\lambda$ is an eigenvalue, then $\lambda >0$: indeed, multiplying \eqref{pbm-vp} by $y$ and integrating by parts, then
$$ \lambda \int _0 ^1 y^2 = \int _0 ^1 x^\alpha y_x ^2 ,$$
which implies first $\lambda \geq 0$, and next that $y=0$ if $\lambda =0$.

%%%%%%%%%%%%

\subsection{The link with the Bessel's equation} \hfill

There is a change a variables that allows one to transform the eigenvalue problem \eqref{pbm-vp} into a
differential Bessel's equation (see in particular Kamke \cite[section  2.162, equation (Ia), p. 440]{Kamke},
and Gueye \cite{Gueye}): 
assume that $\Phi$ is a solution of \eqref{pbm-vp} associated to the eigenvalue $\lambda$; then one easily checks that
the function $\Psi$ defined by
\begin{equation}
\label{eq-lien}
\Phi  (x) =: x^{\frac{1-\alpha}{2}} \Psi \Bigl(\frac{2}{2-\alpha} \sqrt{\lambda} x^{\frac{2-\alpha}{2}} \Bigr)
\end{equation}
is solution of the following boundary problem:
\begin{equation}
\label{pb-bessel}
\begin{cases}
y^2 \Psi ''(y) + y \Psi  '(y) + (y^2 - (\frac{\alpha -1}{2-\alpha}) ^2) \Psi (y) = 0, \quad y\in (0, \frac{2 \sqrt{\lambda} }{2-\alpha}), \\
y^{\frac{1-\alpha}{2-\alpha}} \Psi (y) \to 0 \text{ as } y \to 0  , \\
%(2-\alpha) y^{\frac{1}{2-\alpha}} \Psi  '(y) - (\alpha -1) y^{\frac{\alpha -1}{2-\alpha}} \Psi (y) \to 0 \text{ as } y \to 0 \text{ if } \alpha \in [1,2) , \\
\Psi (\frac{2 \sqrt{\lambda} }{2-\alpha}) = 0 .
\end{cases}
\end{equation}

%%%%%%%%%%%%%%%%%%%%%%%

\subsection{Bessel's equation and Bessel's functions of order $\nu$} \hfill

Bessel's functions of order $\nu$ are solutions of the following differential equation (see \cite[section 3.1, eq. (1), p. 38]{Watson} or \cite[eq (5.1.1), p. 98]{Lebedev}): 
\begin{equation}\label{eq-bessel-ordre-nu}
y^2 \psi''(y) + y \psi'(y) +(y^2-\nu^2) \psi(y)=0, \qquad y\in (0,+\infty).
\end{equation}
The above equation is called {\it Bessel's equation for functions of order $\nu$}. 
Of course the fundamental theory of ordinary differential equations says that the solutions of \eqref{eq-bessel-ordre-nu} generate a vector space $S_\nu$ of dimension 2. Because of \eqref{pb-bessel}, we are interested in solving \eqref{eq-bessel-ordre-nu} when $\nu = \frac{1-\alpha }{2-\alpha}$, hence when $\nu \in (0,\frac{1}{2}]$. When $\nu \notin \Bbb N$, 
\begin{equation}
\label{def-Jnu}
 J_\nu (y)
:= \sum_{m = 0} ^\infty  \frac{(-1)^m}{m! \ \Gamma (m+\nu+1) } \left( \frac{y}{2}\right) ^{2m+\nu}
=\sum_{m = 0} ^\infty  c_{\nu,m} ^+ y^{2m+\nu} 
\end{equation}
and 
\begin{equation}
\label{def-J-nu}
 J_{-\nu} (y)
:= \sum_{m = 0} ^\infty  \frac{(-1)^m}{m! \ \Gamma (m-\nu+1) } \left( \frac{y}{2}\right) ^{2m-\nu}
=\sum_{m = 0} ^\infty  c_{\nu,m} ^- y^{2m-\nu} 
\end{equation}
are well-defined on $\Bbb R^* _+$, and are linearly independent solutions of \eqref{eq-bessel-ordre-nu}. Hence the pair $(J_\nu,J_{-\nu})$ forms a fundamental system of solutions of  \eqref{eq-bessel-ordre-nu},
(see  \cite[section 3.12, eq. (2), p. 43]{Watson}). 

The function $J_\nu$ defined by \eqref{def-Jnu} is the so-called Bessel function of order $\nu$ and  of the first kind. $J_\nu$ has an infinite number of real zeros  which are simple with the possible exception of $x=0$ (\cite[section 15.21, p. 478-479 applied to $C_\nu=J_\nu$]{Watson} or \cite[section 5.13, Theorem 2, p. 127]{Lebedev}). We denote by $(j_{\nu,n})_{n\geq 1}$ the strictly increasing sequence of the positive zeros of $J_{\nu}$:
$$ 0< j_{\nu,1} < j_{\nu,2} < \dots < \ j_{\nu,n} < \dots$$
and we recall that 
$  j_{\nu,n} \to +\infty$ as $n \to +\infty$.

%%%%%%%%%%%%%%%%%%

\subsection{Eigenvalues and eigenfunctions} \hfill

Consider $\alpha \in [0,1)$, and let $\Phi$ be the solution of \eqref{pbm-vp} associated to the eigenvalue $\lambda$. Define $\nu  _\alpha :=  \frac{ 1- \alpha }{2-\alpha} \in (0,\frac{1}{2}]$. 
Hence $\nu  _\alpha \notin \Bbb N$, and Bessel's functions $J_{\nu _\alpha}$ and $J_{-\nu _\alpha}$
are particular solutions of
\begin{equation}
\label{bessel-restr}
y^2 \Psi ''(y) + y \Psi  '(y) + (y^2 - (\frac{\alpha -1}{2-\alpha}) ^2) \Psi (y) = 0 , \quad y \in (0, \frac{2 \sqrt{\lambda} }{2-\alpha}) .
\end{equation}
Since they are also linearly independent, all the solutions of equation \eqref{bessel-restr} are linear combination of $J_{\nu _\alpha}$ and $J_{-\nu _\alpha}$. Hence there exists constants $C_+$ and $C_-$ such that
$$ \forall y\in (0, \frac{2 \sqrt{\lambda} }{2-\alpha}), \quad \Psi (y) = C_+ J_{\nu _\alpha} (y) + C_- J_{-\nu _\alpha} (y) .$$
In particular,
$$ \Phi (x) = C_+  x^{\frac{1-\alpha}{2}} J_{\nu _\alpha} (\frac{2}{2-\alpha} \sqrt{\lambda} x^{\frac{2-\alpha}{2}}) + C_-  x^{\frac{1-\alpha}{2}} J_{-\nu _\alpha} (\frac{2}{2-\alpha} \sqrt{\lambda} x^{\frac{2-\alpha}{2}}).$$
Define
\begin{equation*}
\label{base-phi}
\Phi _+ (x) := x^{\frac{1-\alpha}{2}} J_{\nu _\alpha} (\frac{2}{2-\alpha} \sqrt{\lambda} x^{\frac{2-\alpha}{2}}),
\quad
\Phi _- (x) := x^{\frac{1-\alpha}{2}} J_{-\nu _\alpha} (\frac{2}{2-\alpha} \sqrt{\lambda} x^{\frac{2-\alpha}{2}}).
\end{equation*}
Then, using the series expansion of $J_{\nu _\alpha}$ and $J_{-\nu _\alpha}$, one obtains
\begin{equation*}
%\label{serie-phi}
\Phi _+ (x) = \sum _{m=0} ^\infty \tilde{c} _{\nu _\alpha ,m} ^+  x ^{1-\alpha + (2-\alpha) m},
\quad 
\Phi _- (x) = \sum _{m=0} ^\infty \tilde{c} _{\nu _\alpha ,m} ^- x ^{(2-\alpha) m} 
\end{equation*}
with
\begin{equation*}
%\label{coeffs-phi}
 \tilde{c} _{\nu _\alpha ,m} ^+ := c_{\nu _\alpha ,m} ^+ \Bigl( \frac{2}{2-\alpha} \sqrt{\lambda} \Bigr) ^{2m+\nu _\alpha},
\quad 
\tilde{c} _{\nu _\alpha ,m} ^- := c_{\nu _\alpha ,m} ^- \Bigl( \frac{2}{2-\alpha} \sqrt{\lambda} \Bigr) ^{2m-\nu _\alpha} .
\end{equation*}
Next one easily verifies that $\Phi _+, \Phi _- \in H^1 _\alpha (0,1)$: indeed,
$$ \Phi _+ (x) \sim _{0} \tilde{c} _{\nu _\alpha ,0} ^+  x ^{1-\alpha},
\quad x^{\alpha /2} \Phi _+ ' (x) \sim _{0}  (1-\alpha) \tilde{c} _{\nu _\alpha ,0} ^+  x ^{-\alpha /2} ,$$
$$ \Phi _- (x) \sim _{0} \tilde{c} _{\nu _\alpha ,0} ^- ,
\quad x^{\alpha /2} \Phi _- ' (x) \sim _{0}  (2-\alpha) \tilde{c} _{\nu _\alpha ,1} ^-  x ^{1-\alpha /2} .$$
Hence, given $C_+$ and $C_-$, $\Phi = C_+ \Phi _+ + C_- \Phi _- \in H^1 _\alpha (0,1)$.
But the boundary conditions allow us to obtain information on $C_+$ and $C_-$: since $\Phi (x) \to 0$ as $x\to 0$,
and $\tilde{c} _{\nu _\alpha ,0} ^- \neq 0$, we obtain that $C_- =0$. Hence $\Phi = C_+ \Phi _+$
and $\Psi = C_+ J_{\nu _\alpha}$.
Finally, since $\Psi (\frac{2 \sqrt{\lambda} }{2-\alpha}) = 0$, $C_+ J_{\nu  _\alpha} (\frac{2 \sqrt{\lambda} }{2-\alpha}) = 0$. Since $\Phi$ is an eigenfunction, $\Psi$ is non identically zero. Hence $C_+ \neq 0$, and
$\frac{2 \sqrt{\lambda} }{2-\alpha}$ is a zero of $J_{\nu _\alpha}$: there exists an integer $m\geq 1$ such that
$$ \frac{2 \sqrt{\lambda} }{2-\alpha} = j_{\nu _\alpha, m} .$$
So
$$\lambda = (\frac{2-\alpha}{2})^2  j_{\nu _\alpha, m}^2 = \kappa _\alpha ^2  j_{\nu _\alpha, m}^2 .$$
Therefore, if $\Phi$ is an eigenfunction associated to the eigenvalue $\lambda$, then for some $C_+$ and  $m \in \Bbb N$, $m\geq 1$, we have 
$$ \lambda = \kappa _\alpha ^2  j_{\nu _\alpha, m}^2
\quad \text{ and } \quad \Phi (x) = C_+ x^{\frac{1-\alpha}{2}} J_{\nu _\alpha} ( j_{\nu _\alpha, m} x^{\kappa_\alpha}) .$$

Conversely, one easily verifies that, for all $m\geq 1$ and all $C$
$$  \Phi (x) :=C  x^{\frac{1-\alpha}{2}} J_{\nu _\alpha} ( j_{\nu _\alpha, m} x^{\kappa_\alpha}) $$
is solution of \eqref{pbm-vp}.

Now consider $\Phi _{\alpha, n}$ given by \eqref{fp}.
The family $(\Phi_{\alpha, n})_{n\geq 1}$ forms an orthonormal basis of $L^2(0,1)$: indeed, they are the eigenfunctions of the operator $T_\alpha$:
$$ T_\alpha: L^2(0,1) \to L^2 (0,1), \quad f \mapsto T_\alpha (f) :=u_f $$
where $u_f \in D(A) $ is the solution of the problem $-Au_f = f $, and $T_\alpha$ is self-adjoint  and compact (see Appendix in \cite{fatiha}). 
%The operator $T_\alpha$ sends $L^2(0,1)$ into $H^1 _{\alpha,0} (0,1)$ in a bounded way, and the injection of $H^1 _{\alpha,0} (0,1)$ into $L^2 (0,1)$ is compact (see Appendix in \cite{fatiha}). Hence $T_\alpha$ is compact, and self-adjoint, hence the associated eigenfunctions form a complete system of $L^2(0,1)$. 
The fact that their $L^2$ norm is equal to $1$ comes from a classical identity on Bessel functions, see \cite{Lebedev}, formula (5.14.5), p. 129:
\begin{multline*}
%\label{ortho-bis}
\int _0 ^1 \Phi _{\alpha,n} (x) ^2 \, dx 
= \frac{2\kappa_\alpha}{J_{\nu _\alpha}' (j_{\nu _\alpha ,n} )^2} 
\int _0 ^1  x^{1-\alpha} J_{\nu _\alpha} ( j_{\nu _\alpha, n} x^{\kappa_\alpha}) ^2 \, dx
\\
=\frac{2\kappa_\alpha}{J_{\nu _\alpha}' (j_{\nu _\alpha ,n} )^2} \int _0 ^1 \frac{1}{\kappa_\alpha} y J_{\nu _\alpha} ( j_{\nu _\alpha, n} y) ^2  \, dy 
= \frac{2\kappa_\alpha}{J_{\nu _\alpha}' (j_{\nu _\alpha ,n} )^2}\frac{J_{\nu _\alpha}' (j_{\nu _\alpha ,n} )^2}{2\kappa_\alpha} = 1 .
\end{multline*}

%%%%%%%%%%%%%%%%%%

\subsection{Some bounds on $J_\nu$ and on its zeros} \hfill

\subsubsection{Some bounds on $J_\nu$}
We will use the following bounds from Landau \cite{Landau}:
\begin{equation}
\label{eq-Landau}
\forall \nu >0, \forall x >0, \quad \vert J_\nu (x) \vert \leq \frac{1}{\nu ^{1/3}} 
\quad \text{ and } \quad \vert J_\nu (x) \vert \leq \frac{1}{x ^{1/3}} ,
\end{equation}
and the classical asymptotic development (\cite{Lebedev} p. 122, (5.11.6)):
\begin{multline}
\label{Bessel-DAS}
J_\nu (z) = \Bigl(\frac{2}{\pi z}\Bigr) ^{1/2} \cos (z - \frac{\nu \pi}{2} -\frac{\pi}{4}) (1 + O(\frac{1}{\vert z \vert ^2}))
\\
- \Bigl(\frac{2}{\pi z}\Bigr) ^{1/2} \sin (z - \frac{\nu \pi}{2} -\frac{\pi}{4}) (O(\frac{1}{\vert z \vert })) ,
\end{multline}
valid when $\vert \text{arg } z \vert \leq \pi - \delta$.

%%%%%%%%%%%%%%%%%%

\subsubsection{Some bounds on the zeros of $J_\nu$} \hfill

Using McMahon's formula (see \cite[section 15.53, p. 506]{Watson} applied in the case $\theta=0$ i.e. for $C_\nu=J_\nu$), we can give the following asymptotic expansion of the zeros of $J_\nu$ for any fixed $\nu\geq0$:
\begin{equation}\label{expand-jnun}
j_{\nu,n} = \left( n +\frac{\nu}{2} -\frac{1}{4} \right) \pi  
- \frac{4 \nu^2 -1 }{8  \left( n +\frac{1}{2}\nu -\frac{1}{4} \right) \pi     } 
+ O \left( \frac{1}{n^3} \right) \ \text{ as } n \to +\infty.
 \end{equation}

We will also use the following bounds on the zeros, proved in Lorch and Muldoon \cite{Lorch}:
\begin{equation}
\label{eq-Lorch}
\forall \nu \in [0, \frac{1}{2}], \forall n\geq 1, \quad 
\pi (n + \frac{\nu}{2}-\frac{1}{4}) \leq j_{\nu, n} \leq \pi (n + \frac{\nu}{4}-\frac{1}{8}) 
%, \\
%\forall \nu \geq \frac{1}{2}, \forall n\geq 1, \quad 
%\pi (n + \frac{\nu}{4}-\frac{1}{8}) \leq j_{\nu, n} \leq\pi (n + \frac{\nu}{2}-\frac{1}{4})
 .
\end{equation}

%We also mention the following asymptotic development of the first zero $j_{\nu,1}$  of $J_\nu$ with respect to $\nu$ (\cite[section 15.81, p. 516]{Watson}):
%$$ j_{\nu,1} = \nu + 1,855757 \nu ^{1/3} + O(1) ,$$
%that we will use in section \ref{s6}.

%%%%%%%%%%%%%%%%

%\subsection{Fourier-Bessel Theory} \hfill

%For all $\nu \geq 0$,  Bessel functions satisfy the following orthogonality property 
%(see \cite[eq. (5.14.4) and (5.14.6), p. 129]{Lebedev}):
%\begin{equation}\label{ortho}
%\int_0^1 x J_\nu (j_{\nu,n}x) J_\nu (j_{\nu,m}x) dx 
%= 
%\begin{cases}
%\frac{1}{2}
% [J_{\nu+1} (j_{\nu,n} ) ]^2 & \text{ if } n = m,\\
 %0 & \text{ if } n \neq m.
%\end{cases}
%\end{equation}
%Moreover one can expand arbitrary functions in term of series involving Bessel functions (see \cite[section 18.24, p. 591]{Watson} or 
%\cite[section 5.14, Theorem 3, p. 129]{Lebedev}): 
%\begin{Theorem}
%Let $f(x)$ be a piecewise continuous function in $(0,1)$ of bounded variation in every closed subinterval. Assume that the integral 
%$$\int_0^1 \vert f(x)\vert \sqrt{x} dx$$
%is finite, and consider

%$$ c_m := \frac{2}{J_{\nu +1} (j_{\nu,m})^2} \int _0 ^1 tf(t) J_\nu (j_{\nu,m}t) \, dt . $$
%Then the Fourier-Bessel series 
%$$\sum_{m\geq 1} c_m J_\nu (j_{\nu,m} x) , \qquad 0<x<1$$
%converges to $f(x)$ at every point of continuity of $f$ and to $(f(x+0)+f(x-0))/2$ at every discontinuity point of $f$. 
%\end{Theorem}

%*************************************

%%%%%%%%%%%%%%%%%%%%%%%%%%%%%%%%%
%%%%%%%%%%%%%%%%%%%%%%%%%%%%%%%%%
%%%%%%%%%%%%%%%%%%%%%%%%%%%%%%%%%

%%%%%%%%%%%%%%%%%%%%%%%%%%%%%%%%%%%

\section{Proof of Theorem \ref{thm2}}\label{s7}
 
Let $\alpha \in [0,1)$ be given and consider $T>0$ and $u_0\in L^2(0,1)$. 
Following the ideas of \cite{Lagnese} (in the context of the wave equation), we may reduce the control problem \eqref{eq-u-G} to a moment problem. Then, we will solve this moment problem, using ideas and results of \cite{FR1,FR2}, and of course properties of the eigenvalues and eigenfunctions given in section \ref{s3}.

%%%%%%%%%%%%%

\subsection{Reduction to a moment problem}

In this part, we analyse the problem with formal computations. First, we expand the initial condition $u_0 \in L^2(0,1)$ and the target $u_T \in L^2(0,1)$:
there exists $(\mu_{\alpha ,n}^0)_{n\geq 1}$,  $(\mu_{\alpha ,n}^T)_{n\geq 1} \in \ell ^2(\mathbb N^\star )$ such that
$$ u_0 (x) =\sum _{n\geq 1} \mu_{\alpha ,n} ^0 \Phi_{\alpha ,n} (x),
\quad u_T (x) =\sum _{n\geq 1} \mu_{\alpha ,n}^T \Phi_{\alpha ,n} (x),  \qquad x \in (0,1).
$$
Next we expand the solution $u$ of \eqref{eq-u-G}: 
$$ u(x,t)= \sum_{n\geq 1} \beta_{\alpha ,n} (t) \Phi_{\alpha ,n} (x), \qquad x\in (0,1), \ t \geq 0$$
with
$$\sum_{n\geq 1 } \beta_{\alpha ,n} (t)^2 <+\infty .$$
Therefore the controllability condition $u(\cdot,T) =u_T$ becomes
$$\forall n \geq 1, \qquad \beta_{\alpha ,n} (T)= \mu_{\alpha ,n} ^T.$$

On the other hand, we observe that 
$w_{\alpha ,n} (x,t):= \Phi_{\alpha ,n} (x) e^{\lambda_{\alpha ,n} (t-T)}$ is solution of the adjoint problem:
\begin{equation}\label{pbm-controle3-adjoint}
\begin{cases}
(w_{\alpha ,n})_t +(x^\alpha (w_{\alpha ,n})_x)_x =0 & \qquad x\in(0,1),\ t>0,\\
w_{\alpha ,n}(0,t)=0, \qquad w_{\alpha ,n} (1,t)=0 & \qquad t>0.
\end{cases}
\end{equation}
A combination of \eqref{eq-u-G} and \eqref{pbm-controle3-adjoint} leads to 
\begin{multline*}
0 = \int_0^T\int_0^1 w_{\alpha ,n} (u_t -(x^\alpha u_x)_x) + u ((w_{\alpha ,n})_t+(x^\alpha (w_{\alpha ,n})_x)_x) \\
= \int_0^1 [ w_{\alpha ,n} u ]_0^T dx - \int_0^T [w_{\alpha ,n} x^\alpha u_x ]_0^1 dt  + \int_0^T [u x^\alpha (w_{\alpha ,n})_x ]_0^1 dt\\
= \int_0^1 u(x,T) \Phi_{\alpha ,n}(x)   dx  - \int_0^1 u(x,0) \Phi_{\alpha ,n}(x)  e^{-\lambda_{\alpha ,n} T } dx
- \int_0^T u(0,t)  (x^\alpha (w_{\alpha ,n})_x)  (0,t)  dt\\
= \beta_{\alpha ,n}(T) -  e^{-\lambda_{\alpha ,n} T }  \mu_{\alpha ,n}^0 - \int_0^T  G(t) e^{\lambda_{\alpha ,n} (t-T) } (x^\alpha \Phi_{\alpha ,n} ')(x=0) dt .
\end{multline*}
It follows that
$$ \beta_{\alpha ,n}(T) =   e^{-\lambda_{\alpha ,n} T }  \mu_{\alpha ,n} ^0 + r_{\alpha ,n} \int_0^T  G(t) e^{- \lambda_{\alpha ,n} (T-t) }  dt $$
where we have set
$$r_{\alpha ,n} = (x^\alpha \Phi_{\alpha ,n} ')(x=0).$$

Hence, the controllability condition $u(\cdot,T) =u_T$ implies that
\begin{equation}
\label{moment-bd}
\forall n \geq 1, \quad  r_{\alpha ,n} \int_0^T   G(t) e^{ \lambda_{\alpha ,n} t }  dt 
= -\mu_{\alpha ,n} ^0 + \mu_{\alpha ,n} ^T e^{\lambda _{\alpha ,n} T} . 
\end{equation}
To prove the existence of such a function $G$, it will be necessary to know if $ r_{\alpha ,n} \neq 0$ for all $n$. We prove this property in the following section.

Moreover, since we want a solution of the moment problem that belongs to $H^1(0,T)$, it will be more interesting to see what its derivative has to satisfy. Integrating by parts, we have
$$\int _0 ^T G(t) e^{\lambda _{\alpha,n}t} \, dt
= [ \frac{1}{\lambda _{\alpha,n}} G(t)  e^{\lambda _{\alpha,n}t} ] _0 ^T
- \int _0 ^T \frac{1}{\lambda _{\alpha,n}} G'(t)  e^{\lambda _{\alpha,n}t} \, dt .$$
Hence the derivative $G'$ has to satisfy
\begin{multline}
\label{moment-bd-G'}
 -  \frac{r_{\alpha ,n}}{\lambda _{\alpha,n}} \int _0 ^T G'(t)  e^{\lambda _{\alpha,n}t} \, dt
\\
= -\mu_{\alpha ,n} ^0 + \mu_{\alpha ,n} ^T e^{\lambda _{\alpha ,n} T} 
- \frac{r_{\alpha ,n}}{\lambda _{\alpha,n}}  \Bigl[ G(T)  e^{\lambda _{\alpha,n}T} - G(0) \Bigr] .
\end{multline}
We will provide a solution of this problem that satisfies $G(0)=0=G(T)$.

%%%%%%%%%%%%%%%%%%%%%

\subsection{The generalized derivative of the eigenfunctions at the degeneracy point}

\begin{Lemma}
\label{lemme-neumann-phi}
The eigenfunctions have the following property:
\begin{equation}
\label{neumann-phi}
\forall n\geq 1, \quad x^\alpha \Phi_{\alpha ,n} ' (x)  \to _{x\to 0}
 \frac{(1-\alpha) \sqrt{2\kappa _\alpha }}{2^{\nu_\alpha} \Gamma(\nu_\alpha+1)}
\frac{ (j_{\nu_\alpha,n})^{\nu_\alpha}}{\vert J'_{\nu_\alpha} (j_{\nu_\alpha,n})\vert}   =: r_{\alpha ,n} .
\end{equation}
This generalized derivative at the degeneracy point satisfies
\begin{equation}
\label{posit}
\forall n\geq 1, r_{\alpha ,n} > 0 ,
\end{equation}
and
\begin{equation}
\label{equiv-ralphan}
r_{\alpha ,n} \sim _{n\to \infty}  \rho _\alpha  j_{\nu_\alpha,n}^{\nu_\alpha + 1/2}
\quad \text{ where } \quad \rho _\alpha =  \frac{(1-\alpha) \sqrt{2\kappa _\alpha }}{2^{\nu_\alpha} \Gamma(\nu_\alpha+1)} \frac{\sqrt{\pi}}{\sqrt{2}} .
\end{equation}
\end{Lemma}

{\it Proof of Lemma \ref{lemme-neumann-phi}.} We recall that
$$\Phi_{\alpha ,n}(x)= \varphi_n x^{(1-\alpha)/2} J_{\nu_\alpha} (j_{\nu_\alpha,n} x^{\kappa _\alpha})
\ \text{ with } 
\varphi_n = \frac{\sqrt{2\kappa _\alpha }}{\vert J'_{\nu_\alpha} (j_{\nu_\alpha,n})\vert} . $$
Hence
\begin{multline*}
 x^\alpha \Phi'_{\alpha ,n}(x)=
 \varphi_n \frac{1-\alpha}{2} 
  x^{(\alpha-1)/2} J_{\nu_\alpha}(j_{\nu_\alpha,n} x^{1-\alpha/2})
 \\
+ \varphi_n j_{\nu_\alpha,n} \left(1-\frac{\alpha}{2}	\right) 
 x^{1/2} J'_{\nu_\alpha} (j_{\nu_\alpha,n} x^{1-\alpha/2}) ,
\end{multline*}
and using once again the series expansion of $J_{\nu_\alpha}$ given by \eqref{def-Jnu}, we obtain \eqref{neumann-phi}.

We see directly from the formula that $r_{\alpha ,n} >0$ for all $n \geq 1$. For the limit as $n\to \infty$: we know (see, e.g.,  \cite[remarks on p. 200]{Watson}) that
$$ J_\nu  (y) ^2 + J_{\nu+1} (y) ^2 \sim \frac{2}{\pi y} \quad \text{as } y \to +\infty .$$
Since we always have the relation
$$ J_{\nu+1} (y) = \frac{\nu}{y} J_\nu (y) - J_\nu '(y) ,$$
we obtain that 
\begin{equation}
\label{Delta}
  J'_{\nu _\alpha}( j_{{\nu _\alpha},n} )^2 =  J_{{\nu _\alpha}+1} ( j_{{\nu _\alpha},n} )^2 \sim \frac{2}{\pi j_{{\nu _\alpha},n}} \quad \text{as } n \to \infty .
\end{equation}
Hence
%$$ 
%\frac{ (j_{\nu,n})^\nu}{\vert J'_\nu (j_{\nu,n})\vert} \to + \infty \quad \text{ as } n \to \infty ,$$
%and thus $r_{\alpha ,n} \to + \infty $ as $n\to\infty$. \qed 
$$ r_{\alpha ,n} \sim _{n\to \infty}  \frac{(1-\alpha) \sqrt{2\kappa _\alpha }}{2^{\nu_\alpha} \Gamma(\nu_\alpha+1)} \frac{\sqrt{\pi}}{\sqrt{2}} j_{\nu_\alpha,n}^{\nu_\alpha + 1/2},$$
which is \eqref{equiv-ralphan}. \qed

%%%%%%%%%%%%%%%%%%%%%%%%%%%%%%%%%%

\subsection{Existence and $L^2$-bound for the biorthogonal sequence}\label{s5}

In order to solve the moment problem \eqref{moment-bd}, we will use a sequence $(\sigma_{\alpha, n})_{n\geq 1}$ in $L^2(0,T)$ which is biorthogonal to $(e^{\lambda_{\alpha ,n} t })_{n\geq 1}$, that is 
$$ 
 \int _0 ^T \sigma_{\alpha, n} (t) e^{\lambda_{\alpha ,m} t } \, dt = \delta _{nm}=
\begin{cases} 1 \text{ if } m=n , \\
0 \text{ if } m \neq n .
\end{cases}
$$

The existence of such a sequence follows from general results of Fattorini and Russell \cite{FR1, FR2}. More precisely, we are going prove the following

\begin{Theorem}
\label{thm-biortho1}
Let $(\lambda_{\alpha,n})_{n\geq 1}$ be defined by \eqref{vp}. 
Then there exist positive constants denoted $B_T$ (depending on $T$) and $K$ (independent of $T$), both independent of $\alpha \in [0,1)$, and 
a sequence $(\sigma_{\alpha , n})_{n\geq 1}$ of functions of $ L^2(0,T)$ satisfying the following properties:
\begin{equation}
\label{**2} 
\forall n, m \geq 1, \quad \int_0^T \sigma_{\alpha , n} (t) e^{\lambda_{\alpha,m} t } dt = \delta_{nm} , 
\end{equation}
\begin{equation}
\label{**1} 
\forall n \geq 1, \quad \int_0^T \sigma_{\alpha , n} (t) dt = 0 , 
\end{equation}
and the $L^2$-bounds
\begin{equation}
\label{**3}
\forall n\geq 1, \quad \Vert \sigma_{\alpha , n} \Vert_{L^2(0,T)} \leq B_T e^{K \sqrt{\lambda _{\alpha ,n}}} e^{-\lambda_{\alpha ,n} T} .
\end{equation}
\end{Theorem}

\begin{Remark} {\rm The sequence $(\sigma_{\alpha , n})_{n\geq 1}$ will be the basis that allows us to write a solution $G$ of the moment problem \eqref{moment-bd}. The $L^2$-bounds \eqref{**3} will be useful to ensure the convergence of the associated series giving the control $G$. The orthogonality condition \eqref{**1} is interesting to construct a control $G$ in $H^1 (0,T)$. Finally, the fact that the constants $B_T$ and $K$ are independent of $\alpha \in [0,1)$ will allow us to estimate the controllability cost (Theorem \ref{thm-cost-bd}). }
\end{Remark}

\noindent {\it Proof of Theorem \ref{thm-biortho1}.} 

To have the existence of the biorthogonal family and an estimate of the $L^2$-norm explicit with respect to the degeneracy parameter $\alpha$, we will use results from \cite{FR2}:
\begin{Theorem} (See \cite[Theorems 1.2 and 1.5]{FR2})
\label{thm-unif}
Let $(\lambda _n)_n$ be a sequence such that, for some $\ell >0$,
\begin{equation}
\label{hyp-unif}
\sqrt{\lambda _0} \geq \ell , \quad
\text{ and, for all $n\geq 0$,} \quad \sqrt{\lambda _{n+1}} - \sqrt{\lambda _{n}}\geq \ell .
\end{equation}
Then there is a sequence $(\tilde \sigma _n)_{n\geq 0}$ that is biorthogonal to the family $(e^{-\lambda_n t})_{n \geq 0}$ in $L^2(0,T)$. Moreover, there exist some constants $B(T,\ell)$ and $K(\ell)$ such that
\begin{equation}
\forall n \geq 0, \quad \Vert \tilde \sigma _n \Vert _{L^2(0,T)} \leq B(T,\ell) e^{K(\ell) \sqrt{\lambda _n}}  .
\end{equation}
(The constants $B(T,\ell)$ and $K(\ell)$ are independent of the sequence $(\lambda _n)_n$ as long as assumptions \eqref{hyp-unif} are satisfied.)
\end{Theorem}

Theorem \ref{thm-unif} can be applied directlly to prove \eqref{**2} and \eqref{**3}: when $\alpha \in [0,1)$, $\nu_\alpha \in (0, \frac{1}{2}]$. Thanks to \eqref{eq-Lorch}, we have first
$$\sqrt{\lambda _{\alpha, 1}}= \kappa_\alpha j_{\nu_\alpha , 1} \geq \kappa _\alpha \pi (\frac{3}{4} + \frac{\nu_\alpha}{2}) = \pi [ \frac{3}{8}(2-\alpha) + \frac{1}{4} ( 1-\alpha )] \geq   \frac{3\pi }{8} .$$
Next we turn to the gap between $\sqrt{\lambda _{\alpha, n}}$ and $\sqrt{\lambda _{\alpha, n+1}}$: since
$$ \sqrt{\lambda _{\alpha, n+1}} - \sqrt{\lambda _{\alpha, n}}
= \kappa _\alpha (j_{\alpha, n+1} - j_{\alpha, n}) ,$$
once again using \eqref{eq-Lorch} we have
\begin{multline*}
 \kappa _\alpha (j_{\alpha, n+1} - j_{\alpha, n}) \geq \kappa _\alpha \pi (n+1 + \frac{\nu_\alpha}{2}-\frac{1}{4} - n - \frac{\nu_\alpha}{4} + \frac{1}{8}) 
\\
= \pi \kappa _\alpha (\frac{7}{8} + \frac{\nu_\alpha}{4}) 
= \pi ( \frac{7}{16}(2-\alpha) + \frac{1}{8} ( 1-\alpha ) \geq \frac{7\pi}{16} .
\end{multline*}
Hence there exists $\ell >0$ independent of $\alpha \in [0,1)$ such that
\begin{equation}
\label{unif-vp1}
\forall \alpha \in [0,1), \forall n\geq 1, \quad \sqrt{\lambda _{\alpha, 1}} \geq \ell 
\quad \text{ and } \quad 
 \sqrt{\lambda _{\alpha, n+1}} - \sqrt{\lambda _{\alpha, n}} \geq \ell  .
\end{equation}
Then, applying Theorem \ref{thm-unif}, we conclude that there exists a sequence $(\tilde \sigma _{\alpha , n})_{n\geq 1}$, biorthogonal to the family $(e^{-\lambda_{\alpha , n} t})_{n\geq 1}$ in $L^2(0,T)$, and constants $B_T$ and $K$ such that

$$ \forall \alpha \in [0,1), \forall n \geq 1, 
\quad \Vert \tilde \sigma _{\alpha , n} \Vert _{L^2(0,T)} \leq B_T e^{K \sqrt{\lambda _{\alpha ,n}}} .$$
Now define
$$  \forall n\geq 1, \quad  \sigma_{\alpha ,n}(t) = e^{-\lambda_{\alpha ,n} T} \tilde \sigma_{\alpha ,n} (T-t) .$$
Clearly, $ \sigma_{\alpha ,n} \in L^2(0,T)$ and we see that
\begin{multline*}
\int_0^T \sigma_{\alpha ,n} (t) e^{\lambda_{\alpha ,m} t} dt = e^{- \lambda_{\alpha ,n} T}  \int_0^T \tilde  \sigma_{\alpha ,n} (T-t ) e^{\lambda_{\alpha ,m} t} dt
\\
=  e^{(\lambda_{\alpha ,m}- \lambda_{\alpha ,n}) T}  \int_0^T \tilde  \sigma _{\alpha ,n} (\tau ) e^{- \lambda_{\alpha ,m} \tau} d\tau
 = e^{(\lambda_{\alpha ,m}- \lambda_{\alpha ,n}) T} \delta_{mn}=  \delta_{mn}.
\end{multline*}
Hence \eqref{**2} is satisfied.
Next, we note that
$$ \Vert \sigma_{\alpha ,n} \Vert_{L^2(0,T)}  = e^{-\lambda_{\alpha ,n} T}  \Vert \tilde \sigma_{\alpha ,n} \Vert_{L^2(0,T)} .$$
So \eqref{**3} is satisfied as well. However, this construction does not ensure that \eqref{**1} is satisfied. To prove the existence of a sequence satisfying \eqref{**2}-\eqref{**3}, we slightly modify the previous  construction by adding the artificial eigenvalue $\lambda _{\alpha,0}:=0$. Let
$$ \forall n \geq 0, \lambda ^*_{\alpha ,n} := \lambda _{\alpha ,n} + 1 .$$
Then  $$\sqrt{ \lambda ^* _{\alpha ,0}} \geq 1 $$
and
\begin{multline*} 
\sqrt{\lambda ^*_{\alpha ,n+1}} - \sqrt{\lambda ^* _{\alpha ,n}}
= \frac{\lambda ^* _{\alpha ,n+1} - \lambda ^* _{\alpha ,n}}{\sqrt{\lambda ^* _{\alpha ,n+1}} + \sqrt{\lambda ^* _{\alpha ,n}}} 
\\
= \frac{\lambda _{\alpha ,n+1} - \lambda _{\alpha ,n}}{\sqrt{\lambda ^* _{\alpha ,n+1}} + \sqrt{\lambda ^* _{\alpha ,n}}} 
=  \frac{\sqrt{\lambda _{\alpha ,n+1}} + \sqrt{\lambda _{\alpha ,n}}}{\sqrt{\lambda _{\alpha ,n+1} +1} + \sqrt{\lambda _{\alpha ,n} +1}} (\sqrt{\lambda _{\alpha ,n+1}} - \sqrt{\lambda _{\alpha ,n}}) 
\\
\geq \frac{\sqrt{\lambda _{\alpha ,n+1}} + \sqrt{\lambda _{\alpha ,n}}}{\sqrt{\lambda _{\alpha ,n+1}} + \sqrt{\lambda _{\alpha ,n}} +2 } (\sqrt{\lambda _{\alpha ,n+1}} - \sqrt{\lambda _{\alpha ,n}}) 
\\
\geq \frac{\sqrt{\lambda _{\alpha ,1}}}{\sqrt{\lambda _{\alpha ,1}} +2 } (\sqrt{\lambda _{\alpha ,n+1}} - \sqrt{\lambda _{\alpha ,n}}) .
\end{multline*}
This last quantity is bounded below by a positive constant independent of $\alpha \in [0,1)$. Hence, we can apply Theorem \ref{thm-unif} to the sequence $(\lambda ^*_{\alpha ,n})_{n\geq 0}$, and we get the existence of a sequence $(\tilde \sigma ^* _{\alpha ,n})_{n\geq0}$ that is biorthogonal to the family $(e^{-\lambda ^*_{\alpha , n} t})_{n\geq 0}$ in $L^2(0,T)$: 
\begin{equation}
\label{rajout*}
\forall m, n\geq 0, \quad \int _0 ^T  \tilde \sigma ^* _{\alpha ,n} (t) e^{-\lambda_{\alpha , m} t} e^{-t} \, dt 
= \delta _{nm} ,
\end{equation}
and there exists $B^* _T$ and $K^*$ (both independent of $\alpha \in [0,1)$) such that 
\begin{equation}
\label{rajout**}\forall n \geq 0, \quad \Vert \tilde \sigma _n ^* \Vert _{L^2(0,T)} \leq B_T ^* e^{K^* \sqrt{\lambda ^* _{\alpha, n}}} .
\end{equation}
So the sequence $( \tilde \sigma ^* _{\alpha ,n} (t) e^{-t} )_{n\geq 1}$ is biorthogonal to $(e^{-\lambda_{\alpha , n} t})_{n\geq 1}$ in $L^2(0,T)$, and, applying \eqref{rajout*} to $n\geq 1$ and $m=0$, we see that it satisfies
$$ \forall n \geq 1, \quad \int _0 ^T \tilde \sigma ^* _{\alpha ,n} (t) e^{-t}  \, dt = 0 .$$
Define
$$  \forall n\geq 1, \quad  \sigma_{\alpha ,n}(t) = e^{-\lambda_{\alpha ,n} T} \tilde \sigma ^* _{\alpha ,n} (T-t) e^{-(T-t)}.$$
The sequence $(\sigma_{\alpha ,n})_{n\geq 1}$ satisfies \eqref{**2}-\eqref{**3}: indeed, first
\begin{multline*}
\forall n, m \geq 1, \quad \int_0^T \sigma_{\alpha , n} (t) e^{\lambda_{\alpha ,m} t } dt 
= \int _0 ^T e^{-\lambda_{\alpha ,n} T} \tilde \sigma ^* _{\alpha ,n} (T-t) e^{-(T-t)} e^{\lambda_{\alpha, m} t } dt 
\\
=\int _0^T e^{-\lambda_{\alpha ,n} T} \sigma ^* _{\alpha ,n} (s) e^{-s} e^{\lambda_{\alpha, m} (T-s) } ds 
\\
= e^{-\lambda_{\alpha ,n} T} e^{\lambda_{\alpha, m} T} \int _0 ^T \sigma ^* _{\alpha ,n} (s) e^{-(\lambda_{\alpha, m} +1)s} \, ds 
= \delta_{nm} ,
\end{multline*}
which is \eqref{**2};
next
\begin{multline*}
\forall n \geq 1, \quad \int_0^T \sigma_{\alpha , n} (t) dt =\int _0 ^T e^{-\lambda_{\alpha ,n} T} \tilde \sigma ^* _{\alpha ,n} (T-t) e^{-(T-t)} \, dt
\\
= e^{-\lambda_{\alpha ,n} T} \int _0 ^T \tilde \sigma ^* _{\alpha ,n} (s) e^{-s} \, ds =0   ,
\end{multline*}
which is \eqref{**1};
finally, concerning the $L^2$-bounds we have
\begin{multline*}
\Vert \sigma_{\alpha , n} \Vert_{L^2(0,T)}^2 
= \int _0 ^T e^{-2\lambda_{\alpha ,n} T} \tilde \sigma ^* _{\alpha ,n} (T-t)^2  e^{-2(T-t)} \, dt
\\
= e^{-2\lambda_{\alpha ,n} T}  \int _0 ^T \tilde \sigma ^* _{\alpha ,n} (s)^2  e^{-2s} \, ds
\leq e^{-2\lambda_{\alpha ,n} T} \Vert \sigma ^* _{\alpha , n} \Vert_{L^2(0,T)}^2 .
\end{multline*}
Hence, using \eqref{rajout**}, we get
\begin{multline*}
 \Vert \sigma_{\alpha , n} \Vert_{L^2(0,T)}
\leq e^{-\lambda_{\alpha ,n} T} \Vert \sigma ^* _{\alpha , n} \Vert_{L^2(0,T)}
\leq e^{-\lambda_{\alpha ,n} T} B_T ^* e^{K^* \sqrt{\lambda ^* _{\alpha ,n}}} 
\\
\leq e^{-\lambda_{\alpha ,n} T} B^* _T e^{K^* (1+ \sqrt{\lambda _{\alpha ,n}})} 
= B^* _T e^{K^*} e^{K^* \sqrt{\lambda _{\alpha ,n}}} e^{-\lambda_{\alpha ,n} T},
\end{multline*}
which proves \eqref{**3}. The proof of Theorem \ref{thm-biortho1} is complete. \qed

%%%%%%%%%%%%%%%%%%%%%%%%%%%%%%%%%%%

\subsection{Formal solution of the moment problem}

We have seen that for all $n\geq 1$, $r_{\alpha ,n} \neq 0$ and $r_{\alpha ,n} \to +\infty$ as $n\to \infty$.
Hence, we can consider 
$$ G (t): = \sum _{m\geq 1}\frac{1 }{r_{\alpha ,m}}(-\mu_{\alpha ,m} ^0 + \mu_{\alpha ,m} ^T e^{\lambda _{\alpha ,m} T}) \sigma_{\alpha ,m} (t), $$
where $(\sigma_{\alpha ,m})_{m\geq 1}$ is given by Theorem \ref{thm-biortho1}. 
Formally, $G$ solves the moment problem \eqref{moment-bd}: given $n\geq1$, we have

\begin{multline*}
 r_{\alpha ,n} \int_0^T   G_\alpha (t) e^{ \lambda_{\alpha ,n} t }  dt 
\\
= \sum _{m\geq 1}\frac{1 }{r_{\alpha ,m}}(-\mu_{\alpha ,m} ^0 + \mu_{\alpha ,m} ^T e^{\lambda _{\alpha ,m} T})  r_{\alpha ,n} \int_0^T  \sigma_{\alpha ,m} (t)  e^{ \lambda_{\alpha ,n} t }  dt  
\\
= \sum _{m\geq 1}\frac{1 }{r_{\alpha ,m}}(-\mu_{\alpha ,m} ^0 + \mu_{\alpha ,m} ^T e^{\lambda _{\alpha ,m} T})  r_{\alpha ,n} \delta _{nm} 
\\
= -\mu_{\alpha ,n} ^0 + \mu_{\alpha ,n} ^T e^{\lambda _{\alpha ,n} T}.
\end{multline*}
However, since we want a solution of the moment problem that belongs to $H^1(0,T)$, it will be more interesting to look for a solution of \eqref{moment-bd-G'}.
Consider 
\begin{equation}
\label{g-alpha}
 g_\alpha (t) :=\sum _{m =1}^\infty  \frac{\lambda _{\alpha,m}}{r_{\alpha ,m}} \Bigl(  \mu_{\alpha ,m} ^0 - \mu_{\alpha ,m} ^T e^{\lambda _{\alpha ,m} T} \Bigr)  \sigma _{\alpha,m} (t) , 
\end{equation}
and
\begin{equation}
\label{gd-g-alpha}
G_\alpha (t) := \int _0 ^t g_\alpha (s) \, ds .
\end{equation}
If $g_\alpha \in L^2(0,T)$, then $G_\alpha \in H^1(0,T)$, $G_\alpha(0)=0$, $G_\alpha '(t)=g_\alpha(t)$ a.e., and at least formally $G_\alpha(T)=0$ since all the functions $\sigma _{\alpha, n}$ are of zero mean value. Moreover, 

\begin{multline*}
 -  \frac{r_{\alpha ,n}}{\lambda _{\alpha,n}} \int _0 ^T G_\alpha '(t)  e^{\lambda _{\alpha,n}t} \, dt
=  -  \frac{r_{\alpha ,n}}{\lambda _{\alpha,n}} \int _0 ^T g_\alpha(t)  e^{\lambda _{\alpha,n}t} \, dt
\\
=  -  \frac{r_{\alpha ,n}}{\lambda _{\alpha,n}}  \sum _{m =1}^\infty  \frac{\lambda _{\alpha,m}}{r_{\alpha ,m}} \Bigl(  \mu_{\alpha ,m} ^0 - \mu_{\alpha ,m} ^T e^{\lambda _{\alpha ,m} T} \Bigr)  \int _0 ^T \sigma _{\alpha,m} (t)  e^{\lambda _{\alpha,n}t} \, dt
\\
=  -  \frac{r_{\alpha ,n}}{\lambda _{\alpha,n}}  \sum _{m =1}^\infty  \frac{\lambda _{\alpha,m}}{r_{\alpha ,m}} \Bigl(  \mu_{\alpha ,m} ^0 - \mu_{\alpha ,m} ^T e^{\lambda _{\alpha ,m} T} \Bigr)  \delta _{nm}
\\
=  -  \frac{r_{\alpha ,n}}{\lambda _{\alpha,n}} \frac{\lambda _{\alpha,n}}{r_{\alpha ,n}} \Bigl(  \mu_{\alpha ,n} ^0 - \mu_{\alpha ,n} ^T e^{\lambda _{\alpha ,n} T} \Bigr)
= -\mu_{\alpha ,n} ^0 + \mu_{\alpha ,n} ^T e^{\lambda _{\alpha ,n} T} ,
\end{multline*}
hence \eqref{moment-bd-G'} is satisfied, and at the same time 
\eqref{moment-bd} is satisfied, but now with an $H^1$ function.

%%%%%%%%%%%%%%%%%%%%%%%%%%%%%%%%%%%

\subsection{Rigorous study of the moment problem and of the controllability problem}

\subsubsection{The control $G_\alpha$ belongs to $H^1(0,T)$}

We consider $G_\alpha$ given by \eqref{gd-g-alpha}. We have to check that $G_\alpha$ belongs to $H^1(0,T)$. Let us check that $g_\alpha$ defined by \eqref{g-alpha} belongs to $L^2(0,T)$. First we notice that
$$ g_\alpha (t) = g_\alpha ^0 (t) - g_\alpha ^T (t), $$
where
$$  g_\alpha ^0 (t) :=\sum _{m =1}^\infty  \frac{\lambda _{\alpha,m}}{r_{\alpha ,m}}  \mu_{\alpha ,m} ^0 \sigma _{\alpha,m} (t), $$
and
$$ g_\alpha ^T (t) :=\sum _{m =1}^\infty  \frac{\lambda _{\alpha,m}}{r_{\alpha ,m}}  \mu_{\alpha ,m} ^T e^{\lambda _{\alpha ,m} T}  \sigma _{\alpha,m} (t) . $$
It follows from Lemma \ref{lemme-neumann-phi} that $r_{\alpha,n} \neq 0$ for all $n$, and
$$\frac{\lambda _{\alpha, n}}{ r_{\alpha,n}} \sim _{n\to \infty} \frac{\kappa _\alpha ^2}{\rho _\alpha}  j_{\nu_\alpha,n}^{3/2-\nu_\alpha} .$$
Moreover, thanks to the $L^2$ bounds \eqref{**3}, we see that
$$ \Vert \frac{\lambda _{\alpha,m}}{r_{\alpha ,m}} \mu_{\alpha ,m} ^0 \sigma _{\alpha,m} (t) \Vert _{L^2(0,T)}
\leq \frac{\lambda _{\alpha,m}}{r_{\alpha ,m}} \vert \mu_{\alpha ,m} ^0 \vert B_T e^{K \sqrt{\lambda _{\alpha ,m}}} e^{-\lambda _{\alpha ,m} T}; $$
using \eqref{eq-Lorch}, we obtain that there exists $C_{\alpha,T} \geq 0$ such that
$$ \Vert \frac{\lambda _{\alpha,m}}{r_{\alpha ,m}} \mu_{\alpha ,m} ^0 \sigma _{\alpha,m} (t) \Vert _{L^2(0,T)}
\leq C_{\alpha,T} m^{3/2}  \vert \mu_{\alpha ,m} ^0 \vert e^{K \kappa _\alpha \pi m} e^{-\kappa _\alpha ^2 \pi ^2 (m -\frac{1}{4})^2 T}; $$
when $( \mu_{\alpha ,m} ^0)_m \in \ell ^2 (\Bbb N)$, the series
$$ \sum _{m =1}^\infty m^{3/2}  \vert \mu_{\alpha ,m} ^0 \vert e^{K \kappa _\alpha \pi m} e^{-\kappa _\alpha ^2 \pi ^2 (m -\frac{1}{4})^2 T} $$
is convergent, hence $g_\alpha ^0 \in L^2(0,T)$. Concerning $g_\alpha ^T$: 
we get from \eqref{**3} that
$$ \Vert \frac{\lambda _{\alpha,m}}{r_{\alpha ,m}} \mu_{\alpha ,m} ^T e^{\lambda _{\alpha ,m} T} \sigma _{\alpha,m} (t)\Vert _{L^2 (0,T)} 
\leq \frac{\lambda _{\alpha,m}}{r_{\alpha ,m}} \vert \mu_{\alpha ,m} ^T \vert B_T e^{K \sqrt{\lambda _{\alpha ,m}}} ;$$
using \eqref{eq-Lorch}, we obtain that there exists $C'_{\alpha,T}\geq 0$ such that
$$ \Vert \frac{\lambda _{\alpha,m}}{r_{\alpha ,m}} \mu_{\alpha ,m} ^T e^{\lambda _{\alpha ,m} T} \sigma _{\alpha,m} (t) \Vert _{L^2(0,T)}
\leq C'_{\alpha,T} m^{3/2}  \vert \mu_{\alpha ,m} ^T \vert e^{K \kappa _\alpha \pi m} ; $$
then, when \eqref{decr-target-gene} is satisfied with the constant $K$ given by Theorem \ref{thm-biortho1}, $g_\alpha ^T \in L^2(0,T)$.
Hence, if the initial condition $u_0 \in L^2(0,1)$ and the prescribed target $u_T$
satisfies \eqref{decr-target-gene} with $K$ given by Theorem \ref{thm-biortho1}, then the series defining $g_\alpha$ is convergent in $L^2(0,T)$, and $G_\alpha \in H^1(0,T)$.

%%%%%%%%%%%%%%%%

\subsubsection{The associated solution is driven from the initial state to the prescribed target}

Now, from the definition of the solution $u$ of the boundary control problem \eqref{eq-u-G}, it is natural to consider the problem (see \eqref{eq-v-g})
\begin{equation}
\label{eq-v-g-alpha} 
\begin{cases}
v_t - (x^\alpha v_x)_x = - \frac{p(x)}{p(0)} g_\alpha (t)  , \\
v(0,t)=0=v(1,t), \\
v(x,0)=u_0(x) ,
\end{cases}
\end{equation}
where we recall that 
$$ p(x) = \int _x ^1 \frac{1}{y^\alpha} \, dy .$$
Fix $\varepsilon \in (0,T)$. Then the regularity noted in Remark \ref{rq-reg} allows us to see that
\begin{multline*}
\int _\varepsilon ^T \int _0 ^1 - \frac{p(x)}{p(0)} g_\alpha (t) \Phi _{\alpha ,n} e^{\lambda _{\alpha ,n}t} 
\\
= \int _\varepsilon ^T \int _0 ^1 (v_t - (x^\alpha v_x)_x) \Phi _{\alpha ,n} e^{\lambda _{\alpha ,n}t}  
\\
= [ \int _0 ^1 v \Phi _{\alpha ,n} e^{\lambda _{\alpha ,n}t} ]_\varepsilon ^T
- \int _\varepsilon ^T \int _0 ^1 \lambda _{\alpha ,n} v \Phi _{\alpha ,n} e^{\lambda _{\alpha ,n}t} 
\\
- \int _\varepsilon ^T [ x^\alpha v_x \Phi _{\alpha ,n} e^{\lambda _{\alpha ,n}t}   ]_0 ^1 +  \int _\varepsilon ^T \int _0 ^1 x^\alpha v_x \Phi ' _{\alpha ,n} e^{\lambda _{\alpha ,n}t}
\\
= e^{\lambda _{\alpha ,n}T} \int _0 ^1 v(T) \Phi _{\alpha ,n}
- e^{\lambda _{\alpha ,n}\varepsilon} \int _0 ^1 v(\varepsilon) \Phi _{\alpha ,n} 
- \lambda _{\alpha ,n}  \int _\varepsilon ^T \int _0 ^1 v \Phi _{\alpha ,n} e^{\lambda _{\alpha ,n}t}
\\
+ \int _\varepsilon ^T [ v x^\alpha  \Phi ' _{\alpha ,n} e^{\lambda _{\alpha ,n}t}  ]_0 ^1
- \int _\varepsilon ^T \int _0 ^1 v (x^\alpha \Phi ' _{\alpha ,n})' e^{\lambda _{\alpha ,n}t}
\\
= e^{\lambda _{\alpha ,n}T} \int _0 ^1 v(T) \Phi _{\alpha ,n}
- e^{\lambda _{\alpha ,n}\varepsilon} \int _0 ^1 v(\varepsilon) \Phi _{\alpha ,n}.
\end{multline*}
Letting $\varepsilon \to 0^+$, we obtain  
$$ \int _0 ^T \int _0 ^1 - \frac{p(x)}{p(0)} g_\alpha (t) \Phi _{\alpha ,n} e^{\lambda _{\alpha ,n}t} 
= e^{\lambda _{\alpha ,n}T} \int _0 ^1 v(T) \Phi _{\alpha ,n}
-  \int _0 ^1 u_0 \Phi _{\alpha ,n}  ,$$
hence
\begin{multline*}
e^{\lambda _{\alpha ,n}T} \int _0 ^1 v(T) \Phi _{\alpha ,n} 
\\
= \mu ^0 _{\alpha ,n} + \Bigl( \int _0 ^T g_\alpha (t) e^{\lambda _{\alpha ,n}t} \, dt  \Bigr) \Bigl(  \int _0 ^1 - \frac{p(x)}{p(0)}  \Phi _{\alpha ,n} \, dx \Bigr)
\\
=  \mu ^0 _{\alpha ,n}
-  \frac{\lambda _{\alpha,n}}{r_{\alpha ,n}} \Bigl(  -\mu_{\alpha ,n} ^0 + \mu_{\alpha ,n} ^T e^{\lambda _{\alpha ,n} T} \Bigr)
\Bigl(  \int _0 ^1 - \frac{p(x)}{p(0)}  \Phi _{\alpha ,n} \, dx \Bigr) .
\end{multline*}

\begin{Lemma}
\label{l-coeffs}
The following identity holds:
\begin{equation}
\label{coeffs}
\forall n \geq 1, \quad \int _0 ^1 \frac{p(x)}{p(0)} \Phi _{\alpha,n} (x) \, dx 
= \frac{r_{\alpha,n}}{\lambda_{\alpha,n}} .
\end{equation}
\end{Lemma}

\noindent {\it Proof of Lemma \ref{l-coeffs}.} 
\begin{multline*}
\int _0 ^1 \frac{p(x)}{p(0)} \Phi _{\alpha,n} (x) \, dx
= \frac{1}{\lambda_{\alpha,n}} \int _0 ^1 \frac{p(x)}{p(0)} (-x^\alpha \Phi ' _{\alpha,n} (x))' \, dx
\\
= \frac{1}{\lambda_{\alpha,n}} \Bigl( [\frac{p(x)}{p(0)} (-x^\alpha \Phi ' _{\alpha,n} (x))] _0 ^1 
- \int _0 ^1 \frac{p'(x)}{p(0)} (-x^\alpha \Phi ' _{\alpha,n} (x)) \, dx \Bigr)
\\
= \frac{1}{\lambda_{\alpha,n}} \Bigl( r_{\alpha,n}
- \frac{1}{p(0)} \int _0 ^1  \Phi ' _{\alpha,n} (x) \, dx \Bigr) 
= \frac{r_{\alpha,n}}{\lambda_{\alpha,n}}. \qed
\end{multline*}
Hence 
$$ \forall n\geq 1, \quad e^{\lambda _{\alpha ,n}T} \int _0 ^1 v(T) \Phi _{\alpha ,n} = \mu_{\alpha ,n} ^T e^{\lambda _{\alpha ,n} T} ,$$
which ensures us that $v(T)=u_T$. Then $u(T)=v(T)=u_T$ using \eqref{u:=v,g} and the fact that $G_\alpha (T)=0$. Hence the proof of Theorem \ref{thm2} is complete. \qed

%%%%%%%%%%%%%%%%%%%%%

\section{Structure of the targets:proof of Propositions \ref{prop-reg-targets} and \ref{prop-always-targets}}
\label{sec-struct}

\subsection{An analyticity result: proof of Proposition \ref{prop-reg-targets}}

Consider 
$$ u_T (x) := \sum _{n=1} ^\infty \mu_{\alpha ,n} ^T \Phi _{\alpha ,n} (x) .$$
We recall that
$$ \Phi_{\alpha, n}(x)=  C_{\alpha ,n}
x^{(1-\alpha)/2} J_{\nu _\alpha} (j_{\nu_\alpha,n} x ^{\kappa_\alpha}) ,
\quad \text{ with } \quad C_{\alpha ,n} = \frac{\sqrt{2 \kappa _\alpha }}{\vert J'_{\nu_\alpha} (j_{\nu_\alpha,n} ) \vert}  .$$
Hence
\begin{equation}
\label{express-target}
u_T (x) = \sum _{n=1} ^\infty \mu_{\alpha ,n} ^T  C_{\alpha ,n}
x^{(1-\alpha)/2} J_{\nu _\alpha} (j_{\nu_\alpha,n} x ^{\kappa_\alpha}) .
\end{equation}

Let us study the behavior of $u_T$ near $0$. Using the series expression of $J_{\nu _\alpha}$, we derive from \eqref{express-target} that

$$ u_T (x)  
=  \sum _{n=1} ^\infty \mu_{\alpha ,n} ^T  C_{\alpha ,n}
x^{(1-\alpha)/2}  \Bigl( \sum _{m=0} ^\infty c_{\nu_\alpha , m} ^+ (j_{\nu_\alpha,n} x ^{\kappa_\alpha}) ^{2m+\nu_\alpha} \Bigr).
$$
Formally, exchanging the sums, we obtain
\begin{multline*}  
u_T (x)
= \sum _{m=0} ^\infty c_{\nu_\alpha , m} ^+x^{\kappa_\alpha (2m+\nu_\alpha) + (1-\alpha)/2} \Bigl(\sum _{n=1} ^\infty \mu_{\alpha ,n} ^T  C_{\alpha ,n}
 j_{\nu_\alpha,n} ^{2m+\nu_\alpha} \Bigr)
\\
= \sum _{m=0} ^\infty c_{\nu_\alpha , m} ^+x^{2m \kappa_\alpha + (1-\alpha)} \Bigl(\sum _{n=1} ^\infty \mu_{\alpha ,n} ^T  C_{\alpha ,n}
 j_{\nu_\alpha,n} ^{2m+\nu_\alpha} \Bigr).
\end{multline*}
This is precisely of the form
\begin{equation}
\label{struct-target2}
u_T (x) = x^{1-\alpha} F_\alpha (x^{\kappa_\alpha}),
\end{equation}
with 
\begin{equation}
\label{def-Fholo}
F_\alpha (z):= \sum _{m=0} ^\infty c_{\nu_\alpha , m} ^+  \Bigl(\sum _{n=1} ^\infty \mu_{\alpha ,n} ^T  C_{\alpha ,n}
 j_{\nu_\alpha,n} ^{2m+\nu_\alpha} \Bigr) z^{2m} .
\end{equation}
We will now provide a rigorous proof of the above reasoning. We will need the following

\begin{Lemma}
\label{l-struct2}
If $( \mu_{\alpha ,n} ^T e^{Kn})_n$ remains bounded for some $K>0$, then the function $F_\alpha$ is holomorphic in the disc $\{ z\in \Bbb C, \vert z \vert < \frac{K}{\pi} \}$ and is even.
\end{Lemma}

\noindent {\it Proof of Lemma \ref{l-struct2}.} 
We recall that, by \eqref{Delta}, we have
$$ C_{\alpha ,n} = \frac{\sqrt{2 \kappa _\alpha }}{\vert J'_{\nu_\alpha} (j_{\nu_\alpha,n} ) \vert} \sim _{n\to\infty} \sqrt{ \kappa _\alpha \pi  j_{\nu_\alpha,n} },$$
hence there exists $C^* \geq 0$ independent of $m$ such that
$$ \vert C_{\alpha ,n} j_{\nu_\alpha,n} ^{2m+\nu_\alpha} \vert 
\leq C^*  j_{\nu_\alpha,n} ^{2m+\nu_\alpha +1/2} 
\leq C^* (\pi n )^{2m+1} ,$$
where we used the fact that \eqref{eq-Lorch} implies that
\begin{equation}
\label{maj-j}
\forall \nu \in (0,\frac{1}{2}], \forall n\geq 1,\quad  j_{\nu,n} \leq \pi n.
\end{equation}
 Let us prove the following estimate:

\begin{Lemma}
\label{l-gamma}
There exists some constant $C$ such that, 
\begin{equation}
\label{gamma}
\forall m\in \Bbb N, \forall K>0, \quad \sum_{n=1} ^\infty n^{2m+1} e^{-Kn} \leq C \frac{ (2m+1)! }{K^{2m+2}}.
\end{equation}
\end{Lemma}

\noindent{\it Proof of Lemma \ref{l-gamma}.} The function $x \mapsto x \mapsto x^{2m+1} e^{-Kx}$ is increasing on $[0, (2m+1)/K]$ and decreasing on $[(2m+1)/K, +\infty)$. Hence, 
$$ \text{ if } n+1 \leq \frac{2m+1}{K}, \quad n^{2m+1} e^{-Kn} \leq \int _n ^{n+1}  x^{2m+1} e^{-Kx} \, dx , $$
and
$$ \text{ if } n \geq \frac{2m+1}{K}, \quad (n+1)^{2m+1} e^{-K(n+1)} \leq \int _n ^{n+1}  x^{2m+1} e^{-Kx} \, dx . $$
Denote $m_0$ the integral part of $(2m+1)/K$. Then
\begin{multline*}
\sum_{n=1} ^\infty n^{2m+1} e^{-Kn} 
= \sum_{n=1} ^{m_0} n^{2m+1} e^{-Kn} 
+ \sum_{n=m_0+1} ^\infty n^{2m+1} e^{-Kn} 
\\
= \sum_{n=1} ^{m_0} n^{2m+1} e^{-Kn} 
+ (m_0+1)^{2m+1} e^{-K(m_0+1)}
+ \sum_{n=m_0+1} ^\infty (n+1)^{2m+1} e^{-K(n+1)} 
\\
\leq \sum_{n=1} ^{m_0} \int _n ^{n+1}  x^{2m+1} e^{-Kx} \, dx  
+ (m_0+1)^{2m+1} e^{-K(m_0+1)} +
\sum_{n=m_0+1} ^\infty \int _{n} ^{n+1}  x^{2m+1} e^{-Kx} \, dx
\\
= \int _1 ^{\infty}  x^{2m+1} e^{-Kx} \, dx + (m_0+1)^{2m+1} e^{-K(m_0+1)} 
\\
= \frac{1}{K^{2m+2}} \int _K ^{\infty}  y^{2m+1} e^{-y} \, dy + (m_0+1)^{2m+1} e^{-K(m_0+1)} 
\\
\leq  \frac{1}{K^{2m+2}} \Gamma (2m+2) +  (m_0+1)^{2m+1} e^{-K(m_0+1)},
\end{multline*}
and the Stirling's formula gives \eqref{gamma}. \qed

Now, using \eqref{gamma}, we see that
$$ \sum _{n=1} ^\infty e^{-Kn}  \vert C_{\alpha ,n} \vert j_{\nu_\alpha,n}^{2m+\nu_\alpha} \leq C^* \pi ^{2m+1} C \frac{ (2m+1)! }{K^{2m+2}}. $$ 
On the other hand, 
$$ \vert c_{\nu_\alpha , m} ^+ \vert \leq \frac{1}{m!^2 \, 4^m} .$$
An easy computation shows that the radius of convergence of the series
$$ \sum _{m=1}^\infty \frac{(2m+1)!}{m!^2 \, 4^m} \frac{ \pi ^{2m+1} }{K^{2m+2}} z^m $$
is $\frac{K^2}{\pi^2}$: indeed,  
$$\frac{\frac{(2m+3)!}{(m+1)!^2 \, 4^{m+1}} \frac{ \pi ^{2m+3} }{K^{2m+4}}}{\frac{(2m+1)!}{m!^2 \, 4^m} \frac{ \pi ^{2m+1} }{K^{2m+2}}}= 
\frac{(2m+3)(2m+2)}{4(m+1)^2}\frac{\pi^2}{K^2} \to _{m\to \infty} \frac{\pi^2}{K^2}.$$
Hence if $\vert z \vert ^2 < \frac{K^2}{\pi^2}$, the series defining $F_\alpha$ is convergent, which concludes the proof of Lemma \ref{l-struct2}. \qed

Note that, if $(\mu ^T _{\alpha ,n})_n$ satisfies the assumption of Proposition \ref{prop-reg-targets},
and $x^{\kappa_\alpha} < \frac{K}{\pi}$, then the series 
$$ \sum _{m=0} ^\infty \vert c_{\nu_\alpha , m} ^+ \vert  x^{2m \kappa_\alpha + (1-\alpha)}  \Bigl(\sum _{n=1} ^\infty \mu ^T _{\alpha ,n}  \vert  C_{\alpha ,n} \vert 
 j_{\nu_\alpha,n} ^{2m+\nu_\alpha} \Bigr)
$$ 
is convergent. Therefore our previous argument is justified and \eqref{struct-target2} is valid, with $F_\alpha$ holomorphic in a neighborhood of $0$. We are going to be a little more precise, proving that 
$F_\alpha$ is in fact holomorphic in the horizontal strip $\{z\in \Bbb C, \vert \Im z \vert < \frac{K}{\pi} \}$.

We note that
$$ J_\nu (x) = \sum _{m=0} ^\infty c_{\nu,m} ^+ x^{2m+\nu} = x^\nu \sum _{m=0} ^\infty c_{\nu,m} ^+ x^{2m}
= x^\nu L_\nu (x)$$
where we denote
$$ L_\nu (z) = \sum _{m=0} ^\infty c_{\nu,m} ^+ z^{2m} ,$$
which is holomorphic in $\Bbb C$. 
Hence, coming back to the expression of $u_T$ we have
\begin{multline*}
 u_T (x) = \sum _{n=1} ^\infty \mu_{\alpha ,n} ^T \Phi _{\alpha ,n} (x) 
= \sum _{n=1} ^\infty \mu_{\alpha ,n} ^T C_{\alpha ,n}
x^{(1-\alpha)/2} J_{\nu _\alpha} (j_{\nu_\alpha,n} x ^{\kappa_\alpha})
\\
=  \sum _{n=1} ^\infty \mu_{\alpha ,n} ^T C_{\alpha ,n}
x^{(1-\alpha)/2} (j_{\nu_\alpha,n} x ^{\kappa_\alpha}) ^{\nu _\alpha} L_{\nu_\alpha} (j_{\nu_\alpha,n} x ^{\kappa_\alpha}) .
\end{multline*}
Note that
$$ \frac{1-\alpha}{2} + \kappa_\alpha \nu _\alpha = 1-\alpha .$$
Hence
$$ u_T (x) = x^{1-\alpha} \sum _{n=1} ^\infty \mu_{\alpha ,n} ^T C_{\alpha ,n}
 j_{\nu_\alpha,n} ^{\nu _\alpha} L_{\nu_\alpha} (j _{\nu_\alpha,n} x ^{\kappa_\alpha}) , $$
which gives that
\begin{equation}
\label{struct-target1}
= x^{1-\alpha} \tilde F_\alpha (x ^{\kappa_\alpha}), 
\end{equation}
with
\begin{equation}
\label{comp-target}
\tilde F_\alpha (z) = \sum _{n=1} ^\infty \mu_{\alpha ,n} ^T C_{\alpha ,n}
 j_{\nu_\alpha,n} ^{\nu _\alpha} L_{\nu_\alpha} (j _{\nu_\alpha,n} z)  .
\end{equation}

Let us prove the following 

\begin{Lemma}
\label{l-struct1}
If the sequence $(\mu_{\alpha ,n} ^T e^{Kn})_{n\geq 1}$ is bounded for some $K>0$, then the function $\tilde F_\alpha$ is holomorphic in the horizontal strip $\{ z\in \Bbb C, \vert \Im z \vert <  \frac{K}{\pi} \}$, and even.
\end{Lemma}

\noindent {\it Proof of Lemma \ref{l-struct1}.} 
We derive from the classical asymptotic development recalled in \eqref{Bessel-DAS} that, given $\delta \in (0, \pi)$, there is some $M_\delta$ such that, if the principal argument of $z$ satisfies $\vert \text{arg } z \vert \leq \pi - \delta$, then we have
$$ \vert J_\nu (z) \vert \leq M_\delta \Bigl(\frac{2}{\pi \vert z \vert }\Bigr) ^{1/2} \Bigl( \vert \cos (z - \frac{\nu \pi}{2} -\frac{\pi}{4}) \vert + \vert \sin (z - \frac{\nu \pi}{2} -\frac{\pi}{4}) \vert \Bigr) .$$
But 
$$ \forall z \in \Bbb C, \quad \vert \cos z\vert \leq e^{\vert \Im z\vert }, \quad \vert \sin z\vert \leq e^{\vert \Im z\vert },$$
hence
$$ \vert J_\nu (z) \vert \leq 2 M_\delta  \Bigl(\frac{2}{\pi \vert z\vert }\Bigr) ^{1/2} e^{\vert \Im z \vert } ,$$
and
$$ \vert L_\nu (z) \vert \leq 2M_\delta \frac{1}{\vert z\vert ^{\nu _\alpha} }  \Bigl(\frac{2}{\pi \vert z\vert }\Bigr) ^{1/2} e^{\vert \Im z \vert } .$$
Therefore
\begin{multline*}
\vert \mu_{\alpha ,n} ^T C_{\alpha ,n}
 j_{\nu_\alpha,n} ^{\nu _\alpha} L_{\nu_\alpha} (j _{\nu_\alpha,n} z)   \vert 
\\
\leq 
2 \vert \mu_{\alpha ,n} ^T \vert C_{\alpha ,n}  j_{\nu_\alpha,n} ^{\nu _\alpha}
 M_\delta \frac{1}{ (j _{\nu_\alpha,n} \vert z\vert ) ^{\nu _\alpha} }  \Bigl(\frac{2}{\pi  j _{\nu_\alpha,n} \vert z\vert }\Bigr) ^{1/2} e^{  j _{\nu_\alpha,n} \vert \Im z \vert }.
\end{multline*}
Since $$ C_{\alpha ,n} = \frac{\sqrt{2 \kappa _\alpha }}{\vert J'_{\nu_\alpha} (j_{\nu_\alpha,n} ) \vert} \sim _{n\to\infty} \sqrt{ \kappa _\alpha \pi  j_{\nu_\alpha,n} },$$
by \eqref{maj-j}, we conclude that there exists $M'_\delta$ such that
$$
\vert \mu_{\alpha ,n} ^T C_{\alpha ,n}
 j_{\nu_\alpha,n} ^{\nu _\alpha} L_{\nu_\alpha} (j _{\nu_\alpha,n} z)   \vert  
\leq M' _\delta  \frac{\vert \mu_{\alpha ,n} ^T \vert }{\sqrt{\vert z \vert}}  e^{\pi n \vert \Im z \vert } .$$
Now, if for some $K$ the sequence $(\mu_{\alpha ,n} ^T e^{Kn})_{n\geq 1}$ is bounded,
then the function $\tilde F_\alpha$ is holomorphic in the strip $\{ \vert \Im z \vert <   \frac{K}{\pi} \}$. \qed 

Of course \eqref{struct-target1} and Lemma \ref{l-struct1} completely prove Proposition \ref{prop-reg-targets}. Note that due to analyticity reasons, $F_\alpha$ and $\tilde F_\alpha$ coincide. \qed

%%%%%%%%%%%%%

\subsection{The question of uniformly reachable targets: proof of Proposition \ref{prop-always-targets}}

Part a) of Proposition \ref{prop-always-targets} follows directly from the reachability result given in Theorem \ref{thm2} and Proposition \ref{prop-reg-targets}. The question is then part b): what can be said about targets $u_T$ that satisfy the assumption \ref{decr-target-gene} for all $\alpha \in [0,1)$ ? Using Proposition \ref{prop-reg-targets}, we get that 
$$ u_T (x)= x^{1-\alpha} F_{\alpha} (x^{\kappa_\alpha}) .$$ 
If $u_T$ is nonzero, then consider $p_\alpha$ the first integer such that $F_\alpha ^{(p_\alpha)} (0) \neq 0$. Then, $$ u_T(x)\sim _{x\to 0} \frac{F^{(p_\alpha)} (0)}{p_\alpha !} x^{\kappa_\alpha p_\alpha + 1-\alpha} .$$
Hence the quantity $\kappa_\alpha p_\alpha + 1-\alpha$ has to remain constant on $[0,1)$. This obliges $\alpha \mapsto p_\alpha$ to be locally constant, but this is not sufficient to ensure that $\kappa_\alpha p_\alpha + 1-\alpha$ remains constant. Therefore $u_T$ has to be identically $0$. \qed

%%%%%%%%%%%%%%%%%%%%%

\section{The cost of null controllability: proof of Theorem \ref{thm-cost-bd}}
\label{sec-cost}

\subsection{Upper bounds of Theorem \ref{thm-cost-bd}} 

Fix $\alpha \in [0,1)$, and $u_0 \in L^2(0,1)$. Then we have constructed an admissible control, that drives the solution $u$ of \eqref{eq-u-G} to $0$ in time $T$:
$$ G_\alpha (t) = \int _0 ^t g_\alpha (s) \, ds ,
\quad \text{ where  } \quad g_\alpha (t) = \sum _{m=1} ^\infty \frac{\lambda _{\alpha ,m} }{r_{\alpha ,m}}\mu_{\alpha ,m} ^0 \sigma_{\alpha , m} (t).$$
Hence we will have
$$ C^{H^1}(\alpha, u_0) \leq \Vert G_\alpha \Vert _{H^1} .$$
It remains to estimate the $H^1$ norm of $G_\alpha$.
Since $G_\alpha (0)=0$, it is equivalent to bound the $L^2$ norm of $g_\alpha$, which can be done as follows:
\begin{multline*}
\Vert g_\alpha \Vert _{L^2(0,T)} \leq \sum _{m=1} ^\infty \frac{\lambda _{\alpha ,m} }{r_{\alpha ,m}} \vert \mu_{\alpha ,m} ^0 \vert \Vert \sigma_{\alpha , m} \Vert  _{L^2(0,T)} 
\\
=  \frac{1}{1-\alpha}\frac{2^{\nu_\alpha} \Gamma(\nu_\alpha+1)}{ \sqrt{2\kappa _\alpha }} 
\sum _{m=1} ^\infty \lambda _{\alpha ,m} \vert \mu_{\alpha ,m} ^0 \vert \frac{ \vert J'_{\nu_\alpha}(j_{\nu_\alpha,m}) \vert }{(j_{\nu_\alpha,m})^{\nu_\alpha}} \Vert \sigma _{\alpha,m} \Vert _{L^2(0,T)} 
\\
\leq \frac{1}{1-\alpha}\frac{2^{\nu_\alpha} \Gamma(\nu_\alpha+1)}{ \sqrt{2\kappa _\alpha }} 
\Bigl( \sum _{m=1} ^\infty \vert \mu_{\alpha ,m} ^0 \vert ^2 \Bigr) ^{1/2}
\\
 \Bigl( \sum _{m=1} ^\infty \kappa_\alpha ^4 j_{\nu_\alpha,m} ^{4-2\nu_\alpha} \vert J'_{\nu_\alpha}(j_{\nu_\alpha,m}) \vert ^2  \Vert \sigma _{\alpha,m} \Vert _{L^2(0,T)}^2  \Bigr) ^{1/2} .
\end{multline*}

We are going to estimate the terms that appear in the last series, with respect to the degeneracy parameter $\alpha$:
first, we have already seen that
$$ \vert J'_{\nu_\alpha}(j_{\nu_\alpha,m}) \vert ^2 = \vert J_{\nu_\alpha +1}(j_{\nu_\alpha,m}) \vert ^2 ;$$
the uniform bound \eqref{eq-Landau} from Landau \cite{Landau} gives us that
$$ \vert J'_{\nu_\alpha}(j_{\nu_\alpha,m}) \vert ^2 \leq \frac{1}{(1+\nu _\alpha) ^{2/3}}, $$
hence 
\begin{equation}
\label{serie-estim1}
 \forall \alpha \in [0,1), \forall m \geq 1, \quad \vert J'_{\nu_\alpha}(j_{\nu_\alpha,m}) \vert ^2 \leq 1 .
\end{equation}
Next, we use the bounds \eqref{eq-Lorch} from Lorch-Muldoon \cite{Lorch} to obtain that 
$$ \forall \alpha \in [0,1), \forall m\geq 1, \quad \pi (m- \frac{1}{4}) \leq j_{\nu_\alpha,m} \leq \pi m .$$
Finally, Theorem \ref{thm-biortho1} gives a uniform estimate on $\Vert \sigma _{\alpha,n} \Vert _{L^2(0,T)}$. Hence we obtain that
\begin{multline*}
\sum _{m=1} ^\infty \kappa_\alpha ^4 j_{\nu_\alpha,m} ^{4-2\nu_\alpha} \vert J'_{\nu_\alpha}(j_{\nu_\alpha,m}) \vert ^2  \Vert \sigma _{\alpha,m} \Vert _{L^2(0,T)}^2 
\\
\leq \kappa_\alpha ^4 \pi ^{4-2\nu_\alpha} \sum _{m=1} ^\infty m^{4-2\nu_\alpha} B_T ^2 e^{2K \sqrt{\lambda _{\alpha ,m}}} e^{-\lambda_{\alpha ,m} T} 
\\
\leq B_T ^2 \kappa_\alpha ^4 \pi ^{4-2\nu_\alpha} \sum _{m=1} ^\infty m^{4-2\nu_\alpha} e^{2K \kappa_\alpha \pi m } e^{-\kappa _\alpha ^2 \pi ^2  (m-1/4)^2 T }  .
\end{multline*}
This series is convergent and bounded uniformly with respect to $\alpha \in [0,1)$. Hence 
there is a universal constant $M_2$ such that 
$$ \Vert g_\alpha \Vert _{L^2(0,T)} \leq \frac{M_2}{1-\alpha} \Vert u_0 \Vert _{L^2(0,1)} ,$$
hence by equivalence of norms, we get that
$$ \Vert G_\alpha \Vert _{H^1(0,T)} \leq \frac{M'_2}{1-\alpha} \Vert u_0 \Vert _{L^2(0,1)} ,$$
which implies that the minimum norm control also satisfies this bound, hence 
$$ C^{H^1} (\alpha ,u_0) \leq \frac{M'_2}{1-\alpha} \Vert u_0 \Vert _{L^2(0,1)}  ,$$
and
$$ C^{H^1} _{bd-ctr} (\alpha) \leq \frac{M'_2}{1-\alpha}.$$
This proves the two upper bounds of Theorem \ref{thm-cost-bd}. \qed

%%%%%%%%%%%%%%%%%%%%%%%

\subsection{Lower bounds of  Theorem \ref{thm-cost-bd}}

Now we are interested in lower bounds. One can estimate from below the norm of the biorthogonal family constructed in Theorem \ref{thm-biortho1}, following Hansen \cite{Hansen}. But that would only provide a bound from below of the control given by the moment method, not of the minimum norm control.

So let us consider $u_0 \in L^2(0,1)$, $u_0 \neq 0$, and $G \in \mathcal U ^{ad} (\alpha, u_0)$ an admissible control. We have already seen that
\begin{equation}
\label{eq-mm}
 r_{\alpha ,n} \int_0^T   G(t) e^{ \lambda_{\alpha ,n} t }  dt = -\mu_{\alpha ,n} ^0 .
\end{equation}
We recall that 
\begin{equation}
\label{coeff}
\mu_{\alpha ,n} ^0 = (u_0, \Phi _{\alpha ,n}) = \int _0 ^1 u_0 (x) 
\frac{\sqrt{2 \kappa _\alpha }}{\vert J'_{\nu_\alpha} (j_{\nu_\alpha,n} ) \vert} 
x^{(1-\alpha)/2} J_{\nu_\alpha} (j_{\nu_\alpha,n} x^{\kappa_\alpha}) \, dx .
\end{equation}
We would like to pass to the limit $\alpha \to 1^-$ in this expression. This follow from a continuity argument.
Consider the function
$$ J: (\nu,x) \mapsto J_\nu (x) 
= \sum_{m = 0} ^\infty  \frac{(-1)^m}{m! \ \Gamma (m+\nu+1) } \left( \frac{x}{2}\right) ^{2m+\nu}.$$
Fix $n \geq 1$. The function $J$ is of class $C^1$ on a neighborhood of $(0, j_{0,n})$. Moreover, 
$$ \frac{\partial J}{\partial x} (0, j_{0,n}) = J'_{0} (j_{0,n}) \neq 0 .$$
Then the implicit function theorem says that there exists a neighborhood 
of $(0, j_{0,n})$ and a function $\psi _n: \mathcal{V} (0) \to \mathcal{V} ( j_{0,n})$ 
such that $J(\nu, \psi _n(\nu))=0$, and 
$$ \begin{cases} J(\nu,x)=0 , \\ (\nu , x) \in \mathcal V (0, j_{0,n})  \end{cases} \text{ is equivalent to } 
 \begin{cases} x=\psi _n (\nu) , \\ \nu  \in \mathcal V (0) . \end{cases}$$
Hence for $\nu >0$ small enough, $J_\nu$ has a zero $\psi_1(\nu)$ close to $j_{0,1}$, a zero $\psi_2(\nu)$ close to $j_{0,2}$, $\cdots$, a zero $\psi_n(\nu)$ close to $j_{0,n}$. But we already know that the zeroes of $J_\nu$ are $j_{\nu,1}, j_{\nu,2}  , \cdots, j_{\nu,n}, \cdots$. 
The bound \eqref{eq-Lorch} says that $j_{\nu,n} \in [\pi(n+\frac{\nu}{2}-\frac{1}{4}), \pi(n+\frac{\nu}{2}-\frac{1}{8})]$, hence there is no choice: for $\nu >0$ small enough, $\psi_1(\nu)=j_{0,1}$, $\cdots$, $\psi_n(\nu)=j_{0,n}$,
and $\nu \mapsto j_{\nu,1}$, $\cdots$, 
$\nu \mapsto j_{\nu,n}$ are continuous in a neighborhood of $0$.

This allows us to pass to the limit in \eqref{coeff}, and we obtain 
$$  \mu_{\alpha ,n} ^0 \to _{\alpha \to 1^-} 
\int _0 ^1 u_0 (x) 
\frac{1}{\vert J'_{0} (j_{0,n} ) \vert} 
 J_{0} (j_{0,n} x^{1/2}) \, dx = (u_0, \Phi _{1 ,n}) ,$$
where we have set
\begin{equation}
\label{Phi1}
\forall n\geq 1, \qquad 
 \Phi_{1, n}(x) :=  \frac{1}{\vert J'_0 (j_{0,n} ) \vert}  J_{0} (j_{0,n} \sqrt{x}), \qquad x\in (0,1) ,
\end{equation}
obtained taking $\alpha =1$ in \eqref{fp}. We note that, as for $\alpha \in [0,1)$, we have the

\begin{Lemma}
\label{thm-gueye-ext}
The family $(\Phi_{1, n})_{n\geq 1}$ forms an orthonormal basis of $L^2(0,1)$.
\end{Lemma}

\noindent {\it Proof of Lemma \ref{thm-gueye-ext}.}
First we check that the family $(\Phi_{1, n})_{n\geq 1}$ forms an orthonormal family of $L^2(0,1)$.
First we consider $n\neq m$. By change of variables, we have
\begin{multline*}
 \int _0 ^1  \Phi_{1, n}(x)  \Phi_{1, m}(x) \, dx 
= \frac{1}{\vert J'_0 (j_{0,n} ) \vert}  \frac{1}{\vert J'_0 (j_{0,m} ) \vert} \int _0 ^1  J_{0} (j_{0,n} \sqrt{x}) J_{0} (j_{0,m} \sqrt{x}) \, dx
\\
= \frac{1}{\vert J'_0 (j_{0,n} ) \vert}  \frac{1}{\vert J'_0 (j_{0,m} ) \vert} \int _0 ^1  2y J_{0} (j_{0,n} y) J_{0} (j_{0,m} y) \, dy .
\end{multline*}
Then we use \cite{Lebedev}, formula (5.14.3) p. 128
to get
$$ \int _0 ^1  2y J_{0} (j_{0,n} y) J_{0} (j_{0,m} y) \, dy 
= \frac{j_{0,m}J_{0} (j_{0,n} )J' _{0} (j_{0,m})-j_{0,n}J_{0} (j_{0,m} )J' _{0} (j_{0,n})}{n^2-m^2}
= 0 .$$
Hence if $n\neq m$, we have
$$ \int _0 ^1  \Phi_{1, n}(x)  \Phi_{1, m}(x) \, dx = 0.$$
When $n=m$, we use \cite{Lebedev}, formula (5.14.5) p. 129:
\begin{multline*}
 \int _0 ^1  \Phi_{1, n}(x) ^2 \, dx 
= \frac{1}{ J'_0 (j_{0,n} ) ^2}  \int _0 ^1  2y J_{0} (j_{0,n} y) ^2  \, dx
\\
= \frac{1}{ J'_0 (j_{0,n} ) ^2} ( J'_0 (j_{0,n} ) ^2 + J_0 (j_{0,n} ) ^2 ) = 1 .
\end{multline*}
At last we check that the family $(\Phi_{1, n})_{n\geq 1}$ generates $L^2(0,1)$. Consider a smooth function $\psi$ compactly supported in $(0,1)$. We would like to prove that the series 
$$ \sum _{n=1} ^\infty (\psi,  \Phi_{1, n})  \Phi_{1, n} (x) $$
converges to $\psi$ in $L^2(0,1)$. We note that
$$
\sum _{n=1} ^\infty (\psi,  \Phi_{1, n})  \Phi_{1, n} (x) 
\\
= \sum _{n=1} ^\infty \Bigl( \frac{2}{J'_0 (j_{0,n} ) ^2} \int _0 ^1 y \psi (y^2) J_{0} (j_{0,n} y) \, dy \Bigr) J_0 (j_{0,n} \sqrt{x}) .
$$
It follows from \cite{Watson} chapter XVIII, sections 18.24-18.26, 18.53 that the series
$$ \sqrt{x} \sum _{n=1} ^\infty \Bigl( \frac{2}{J'_0 (j_{0,n} ) ^2} \int _0 ^1 y \psi (y^2) J_{0} (j_{0,n} y) \, dy \Bigr) J_0 (j_{0,n} x) $$
converges uniformly on $[0,1]$ to $\sqrt{x} \psi (x^2)$. Hence
$$\sup _{x\in[0,1]} \vert x^{1/4} \psi (x) - x^{1/4} \sum _{n=1} ^N \Bigl( \frac{2}{J'_0 (j_{0,n} ) ^2} \int _0 ^1 y \psi (y^2) J_{0} (j_{0,n} y) \, dy \Bigr) J_0 (j_{0,n} \sqrt{x}) \vert $$
goes to $0$ as $N\to \infty$.
Therefore
$$\sup _{x\in[0,1]} \vert x^{1/4} \psi (x) - x^{1/4} \sum _{n=1} ^N (\psi,  \Phi_{1, n})\Phi_{1, n}(x) \vert \to _{N\to \infty} 0 .$$
But
\begin{multline*}
\int _0 ^1 \vert \psi (x) - \sum _{n=1} ^N (\psi,  \Phi_{1, n})  \Phi_{1, n} (x) \vert ^2 \, dx 
\\ 
= \int _0 ^1 \vert x^{1/4} \psi (x) - x^{1/4} \sum _{n=1} ^\infty (\psi,  \Phi_{1, n})(x) \vert ^2 \frac{dx}{\sqrt{x}} . 
\end{multline*}
Since $x^{-1/2} \in L^1(0,1)$, we have
$$
\int _0 ^1 \vert \psi (x) - \sum _{n=1} ^N (\psi,  \Phi_{1, n})  \Phi_{1, n} (x) \vert ^2 \, dx
 \to _{N\to \infty} 0 .
$$
This implies that every smooth compactly supported function $\psi$ is the limit in $L^2(0,1)$ of 
linear combination of the $\Phi _{1,n}$. Since these functions are dense in $L^2(0,1)$, the family $(\Phi_{1, n})_{n\geq 1}$ is an an orthonormal basis of $L^2(0,1)$.
\qed

\begin{Remark} {\rm 
In fact it can be proved that the functions $\Phi_{1, n}$ are the eigenfunctions of the problem
\begin{equation}\label{pbm-vp-N}
\begin{cases}
 - (x y'(x))' =\lambda y(x) & \qquad x\in (0,1),\\
(xy') (0)=0 ,\\ 
y(1)=0 , & 
\end{cases}
\end{equation}
associated to the eigenvalues $\lambda_{1, n} =  \kappa _1 ^2 j_{0 ,n}^2$. We will prove and use this (stronger) property in \cite{CMV-cost-loc}.}
\end{Remark}

Now, we note that since $(\Phi _{1 ,n})_{n\geq1}$ forms a orthonormal basis of $L^2(0,1)$, there is $n$ such that
$$ (u_0, \Phi _{1 ,n}) \neq 0 .$$
Hence, there exists $m_0 (u_0)>0$ and $n_0$ such that, for $\alpha$ sufficiently close to $1^-$, we have
$$ \vert \mu^0 _{\alpha ,n_0} \vert \geq m_0(u_0) .$$
Now, we deduce from \eqref{eq-mm} and the Cauchy-Schwarz inequality that
$$\frac{\vert \mu_{\alpha ,n_0} ^0 \vert }{r_{\alpha ,n_0}} \leq \Bigl( \int _0 ^T e^{2\lambda _{\alpha ,n_0}t} \, dt \Bigr)^{1/2} 
\Bigl( \int _0 ^T G(t)^2 \, dt \Bigr)^{1/2} ,$$
hence
$$ \frac{m_0(u_0) }{r_{\alpha ,n_0}} \leq 
\frac{e^{2\lambda _{\alpha ,n_0} T}-1}{2\lambda _{\alpha ,n_0}} \Vert G \Vert _{L^2(0,T)} .$$
We thus obtain a bound from below for the admissible control:
$$ \Vert G \Vert _{L^2(0,T)}
\geq \frac{m_0(u_0)}{1-\alpha} \frac{2\lambda _{\alpha ,n_0}}{e^{2\lambda _{\alpha ,n_0} T}-1}
 \frac{2^{\nu_\alpha} \Gamma(\nu_\alpha+1)}{ \sqrt{2\kappa _\alpha }}
\frac{\vert J'_{\nu_\alpha} (j_{\nu_\alpha,n_0})\vert}{ (j_{\nu_\alpha,n_0})^{\nu_\alpha}}.$$
Now we conclude noting that
\begin{multline*}
 \vert J'_{\nu_\alpha} (j_{\nu_\alpha,n_0})\vert
= \vert J_{\nu_\alpha +1} (j_{\nu_\alpha,n_0})\vert
= \vert J(\nu_\alpha +1, j_{\nu_\alpha,n_0}) \vert
\\
\to _{\alpha \to 1^-} \vert J(1, j_{0,n_0}) \vert = J_1( j_{0,n_0}) = J_0 ' ( j_{0,n_0}) \neq 0 ,
\end{multline*}
hence $\vert J'_{\nu_\alpha} (j_{\nu_\alpha,n_0})\vert$ is bounded from below by a positive constant. Since $j_{\nu_\alpha,n_0}$ is bounded from below and from above by constants depending on $n_0$ but uniform with respect to $\alpha \in [0,1)$, 
there is $m_1(u_0) >0$, depending on $n_0$ (hence on $u_0$), but independent of $\alpha \in [0,1)$ such that, for all $G \in \mathcal U _{ad} (\alpha,u_0)$ we have
$$ \Vert G \Vert _{L^2(0,T)} \geq \frac{m_1(u_0)}{1-\alpha} .$$
Of course, the $H^1$ norm of $G$ will satisfy the same lower estimate, which concludes the proof of
the lower bound in \eqref{eq-cost-fr-u_0}. Then the lower bound in \eqref{eq-cost-fr} follows and the proof of Theorem  \ref{thm-cost-bd} is complete. \qed

%%%%%%%%%%%%%%%%%%%%%%%%%%%%%%%%%%%%%%%%
%%%%%%%%%%%%%%%%%%%%%%%%%%%%%%%%%%%%%%%%

\end{document}